\documentclass[11pt]{amsart}

\usepackage{amsmath,amssymb,amsthm}
\usepackage[inner=1in,outer=1in,top=1in,bottom=1in]{geometry} 
\usepackage{amsthm,amsmath,amssymb,amscd,graphics,enumerate,stmaryrd,xspace,verbatim,epic,eepic,url}
\usepackage{mathtools}
\usepackage{listings} 
\usepackage{caption} 
\usepackage{tikz-cd} 
\usepackage[colorlinks=true]{hyperref} 
\usepackage{array}   
\usepackage{longtable} 

\newcommand{\sheafhom}{\mathscr{H}\kern -.5pt om}

\newcommand{\OY}{\mathcal{O}_Y}

\newcommand{\into}{\hookrightarrow}

\newcommand{\bu}{\textbullet \ }

\newcommand{\Spec}{\operatorname{Spec}}
\newcommand{\Proj}{\operatorname{Proj}}
\newcommand{\cProj}{{{\mathcal P}roj}}
\newcommand{\Hom}{\operatorname{Hom}}

\newcommand{\inv}{\operatorname{inv}}

\newcommand{\ord}{\operatorname{ord}}
\newcommand{\Der}{\operatorname{Der}}

\newcommand{\maxord}{\operatorname{max-ord}}
\newcommand{\lcm}{\operatorname{lcm}}
\newcommand{\maxinv}{\operatorname{maxinv}}
\newcommand{\ZR}{{\mathbf{ZR}}}
\newcommand{\bmu}{{\boldsymbol{\mu}}}
\newcommand{\Diff}{\operatorname{Diff}}
\newcommand{\sDiff}{\mathcal{D}\textit{iff}}
\newcommand{\lexinv}{\operatorname{inv}}
\newcommand{\binv}{\operatorname{\text{$b-$}inv}}
\newcommand{\btoa}{\operatorname{\text{$b-to-a$}}}

\let\emptyset\varnothing

\newcommand{\MC}{\operatorname{MC}}
\newcommand{\con}{\operatorname{CONTENDERS}}
\newcommand{\surv}{\operatorname{SURVIVORS}}
\newcommand{\win}{\operatorname{WINNERS}}
\newcommand{\lose}{\operatorname{LOSERS}}
\newcommand{\ROW}{\operatorname{ROW}}
\newcommand{\COL}{\operatorname{COL}}

\newcommand{\bmax}{\operatorname{b}_{\operatorname{max}}}
\newcommand{\bcurr}{\operatorname{b}_{\operatorname{curr}}}
\newcommand{\maxbinv}{\operatorname{maxbinv}}
\newcommand{\Res}{\operatorname{RESOLUTION}}

\newtheorem{definition}[subsubsection]{Definition}

\newtheorem{proposition}[subsubsection]{Proposition}
\newtheorem{corollary}[subsubsection]{Corollary}
\newtheorem{lemma}[subsubsection]{Lemma}

\makeatletter \def\@dotsep{4.5} \makeatother 

\renewcommand*\thelstnumber{\arabic{lstnumber}:}
\DeclareCaptionFormat{listing}{\rule{\dimexpr\textwidth\relax}{1pt}\par#1#2#3}

\let\origthelstnumber\thelstnumber
\makeatletter
\newcommand*\Suppressnumber{%
  \lst@AddToHook{OnNewLine}{%
    \let\thelstnumber\relax%
     \advance\c@lstnumber-\@ne\relax%
    }%
}
\newcommand*\Reactivatenumber{%
  \lst@AddToHook{OnNewLine}{%
   \let\thelstnumber\origthelstnumber%
   \advance\c@lstnumber\@ne\relax}%
}

\begin{document}
\title{Algorithmic resolution via weighted blowings up}

\author[Jonghyun Lee]{Jonghyun Lee}
\address{Department of Mathematics, University of Michigan,
Ann Arbor, MI, 48109, U.S.A.}
\email{nuyhgnoj@umich.edu}


\date{\today}

\begin{abstract}
In this paper we describe a computer implementation of Abramovich, Temkin, and W{\l}odarczyk's algorithm for resolving singularities in characteristic zero. Their ``weighted resolution" algorithm proceeds by repeatedly blowing up along centers that are independent of the history of the past blowing ups, distinguishing weighted resolution from previous resolution algorithms, which all rely on history. We compare our implementation of weighted resolution with that of Villamayor's resolution algorithm, experimentally verifying that weighted resolution is remarkably efficient. 
\end{abstract}
\maketitle
\setcounter{tocdepth}{1}

\section{Introduction} 
\label{sec: intro}

In its simplest form, a resolution of singularities of a variety $X$ is a proper birational morphism $Y \to X$ from some smooth variety $Y$. Determining the existence and construction of resolution of singularities is one of the central problems in algebraic geometry. In his landmark paper \cite{Hir64}, Hironaka proved that every variety over a field of characteristic zero admits a resolution of singularities. Through the works of Villamayor, Bravo, Encinas, Bierstone, Milman, and W{\l}odarczyk, \cite{Vil89,Vil92,BM97,EV98,BV01,BEV05,Wlo05}, Hironaka's intricate, non-constructive proof has been vastly simplified to constructive resolution algorithms. But these further developments still rely on Hironaka's original argument of repeatedly blowing up along centers that depend on the cumulative history of the past blowing ups. History forces resolution algorithms to make more and more computations when calculating the next center. Unfortunately, history is unavoidable, that is, as long as we stay in the realm of varieties. 

By moving to a bigger realm which subsumes varieties, the realm of Deligne-Mumford stacks, we can avoid history. Replacing blowing ups by stack-theoretic weighted blowing ups, Abramovich, Temkin, and W{\l}odarczyk discovered a constructive resolution algorithm that drops history and assigns a canonical center based on the input variety alone \cite[Theorem 1.1.1]{ATW20}!\footnote{A similar result was also discovered independently by McQuillan \cite{McQ19}.} Their ``weighted resolution" algorithm is far simpler than previous resolution algorithms, and is also remarkably efficient. 

With Abramovich and Fr{\"u}bis-Kr{\"u}ger, the author implemented the weighted resolution algorithm in the computer algebra system {\sc Singular} \cite{DGPS19}, which will appear as the {\sc Singular} library \texttt{resweighted.lib} \cite{AFL20}. The aim of this paper is to describe this computer implementation. In Section \ref{sec: resweighted}, we provide the algorithm \nameref{weightedResolution}, which is the weighted resolution algorithm. At the end of the paper in Section \ref{sec: compare}, we will compare the weighted resolution algorithm with Villamayor's resolution algorithm, which is implemented as the {\sc Singular} library \texttt{resolve.lib} \cite{FP19}, demonstrating by computer experiments that weighted resolution is indeed remarkably efficient. 

Developing computer implementations of constructive resolution algorithms is highly nontrivial work. Despite the word ``constructive" in the name, constructive resolution algorithms\footnote{See \cite{EV98} for the definition of ``algorithm of resolution," which is what we call constructive resolution algorithm. Also, see \cite[page 630]{Vil92}.} rely on non-constructive existence results which are not useful for computer implementations, since computers need to be told explicitly what to do.\footnote{For example, passing to formal completions or taking an open cover with sufficiently small open sets are common procedures that are not constructed in constructive resolution algorithms.} The main challenges of implementing the weighted resolution algorithm reside in creating explicit constructions for finding the blowing up centers and performing the weighted blowings up.

In weighted resolution, as well as in previous resolution algorithms, determining the blowing up center requires computing maximal contact hypersurfaces. Theoretically, maximal contact hypersurfaces always exist in sufficiently small neighborhoods of singularities, but the explicit computation of these neighborhoods and maximal contact hypersurfaces is one of the main challenges of implementing the weighted resolution algorithm. In Section \ref{sec: maximal contact}, we provide algorithms \nameref{maximalContact} and \nameref{liftMaximalContact} for computing maximal contact hypersurfaces, upgrading the procedure for maximal contact hypersurfaces given in \cite[\texttt{CoverCenter}]{FP19}.\footnote{\texttt{CoverCenter} is a static procedure that can only be viewed in the source code of the {\sc Singular} library \texttt{resolve.lib} \cite{FP19}.} Another approach for obtaining maximal contact hypersurfaces can be found in \cite{BS00a,BS00b}.

After constructing the center, we need to compute the weighted blowing up along that center. In general, the weighted blowings up of a variety is a Deligne-Mumford stack, which is locally the quotient of some variety by a finite group. Thus we can cover the weighted blowings up by affine charts equipped with finite group actions. Local descriptions of these charts and group actions are provided in \cite[Section 3.4]{ATW20}. However, these local descriptions are non-constructive, making the explicit computation of these charts and group actions another one of the main challenges of implementing weighted resolution. In Section \ref{sec: blowup}, we provide algorithm \nameref{weightedBlowUp}, which explicitly computes these charts and group actions.

\subsection{The weighted resolution algorithm}
We define the weighted resolution algorithm to be the procedure for embedded resolution of singularities given by \cite[Theorem 1.1.1]{ATW20}. See Algorithm \nameref{weightedResolution} in Section \ref{sec: resweighted} for the description of the implementation of the weighted resolution algorithm. Each algorithm written in this paper is a procedure (with the same name) in the {\sc Singular} library \texttt{resweighted.lib}. We emphasize that the weighted resolution algorithm would not exist without the prior works of Villamayor, Bierstone, Milman, Hironaka, and many others. Moreover, the implementation of the weighted resolution algorithm in \cite{AFL20} would not exist without the prior computer implementations by Bodn{\'a}r and Schicho in \cite{BS00a,BS00b} and Fr\"{u}hbis-Kr\"{u}ger and Pfister in \cite{FP19}. 


\subsection{Overview of this paper}
In Section \ref{sec: ZR}, we run through the formalisms of Zariski-Riemann spaces, Rees algebras, and valuative $\mathbb{Q}$-ideals. These formalisms are needed in order to precisely define the centers used in weighted resolution (which we do in the same section), allowing us in Section \ref{sec: lexinv} to describe the output of  Algorithm \nameref{prepareCenter}, which is the data of a center, with complete unambiguity. Furthermore, these formalisms are needed again in Section \ref{sec: blowup}, where they supply the proofs that show that Algorithm \nameref{weightedBlowUp} correctly outputs the charts of weighted blowing ups.

Sections \ref{sec: derivatives}, \ref{sec: maximal contact}, and \ref{sec: lexinv} provide the algorithms used to find the center associated to an ideal on a smooth variety. Section \ref{sec: derivatives} provides Algorithms \nameref{Diff}, which computes the first derivative $\mathcal{D}^{\leq 1}_{Y/k} \mathcal{I}$ of a coherent ideal $\mathcal{I}$ on a smooth variety $Y$, and \nameref{maximalOrderOfVanishing}, which computes the maximal order of vanishing of a coherent ideal on a smooth variety; Section \ref{sec: maximal contact} provides Algorithms \nameref{maximalContact} and \nameref{liftMaximalContact}, which computes maximal contact hypersurfaces and their lifts, respectively; and Section \ref{sec: lexinv} provides Algorithm \nameref{prepareCenter}, which computes the center associated to an ideal on a smooth variety using the previous algorithms. 

Section \ref{sec: blowup} provides Algorithm \nameref{weightedBlowUp}, which explicitly computes the charts of weighted blowing ups. Section \ref{sec: resweighted} provides the culmination of all the algorithms in this paper, Algorithm \nameref{weightedResolution}, which explicitly computes the weighted resolution of a singular variety embedded in a smooth variety. 

In Section \ref{sec: compare}, we compare our computer implementation \texttt{resweighted.lib} \cite{AFL20} of weighted resolution with the {\sc Singular} library \texttt{resolve.lib} \cite{FP19}, which implements Villamayor's resolution algorithm. By computer experiments, we demonstrate that weighted resolution is significantly more efficient than Villamayor's algorithm. 

There are two appendices in this paper. Appendix \ref{sec: saturation} establishes the correspondence between Zariski closure and the ideal-theoretic operation of saturation. We include this as an appendix due to the lack of scheme-theoretic treatment in the literature. Also, Algorithm \nameref{Diff} relies on results proven in Appendix \ref{sec: saturation}. Appendix \ref{sec: helper} contains Algorithms \nameref{isSmoothHypersurface}, which determines whether an equation cuts a smooth hypersurface out of a smooth variety, and \nameref{orthogonal idempotents}, which computes the orthogonal idempotents on a smooth variety. Though these two algorithms are not centrally involved in weighted resolution, we include them as an appendix because they are needed for the main algorithms in this paper.

\subsection{Acknowledgements}
I am immensely grateful to Dan Abramovich and Anne Fr\"{u}hbis-Kr\"{u}ger for patiently advising me through developing a computer implementation of the weighted resolution algorithm as well as the writing of this paper. 

\subsection{Conventions}
Let the natural numbers be $\mathbb{N}=\{0,1,2,\dots\}$, and let $\mathbb{Q}_{>0}$ denote the positive rational numbers. Since everybody seems to have their own definition of an algebraic variety, we will state which definition we use. We follow \cite[Definition 2.3.47]{Liu06}, defining a variety over a field $k$ to be a scheme of finite type over $k$.

{\tableofcontents}

\renewcommand\lstlistlistingname{List of Algorithms}
{\lstlistoflistings}

\pagebreak

\section{Zariski-Riemann spaces, Rees algebras, and centers}
\label{sec: ZR}
Let $Y$ be a separated, integral variety over $k$, with function field $K$. We introduce the formalism of Zariski-Riemann spaces, and Rees algebras, valuative $\mathbb{Q}$-ideals, in order to precisely define the centers involved in weighted resolution. We follow \cite{Que20}. 

\subsection{Zariski-Riemann spaces}
 The Zariski-Riemann space $\ZR(Y)$ of $Y$ is a topological space whose points are valuation rings $R_v$ of $K$ which contain $k$ and possess a center on $Y$. $\ZR(Y)$ is quasicompact and irreducible. We will refer to points of $\ZR(Y)$ by valuations $v$ on $K$ rather than their corresponding valuation rings $R_v$. Let $\pi_Y: \ZR(Y) \to Y$ be the continuous map sending a valuation $v$ to its unique center $\pi_Y(v)$ on $Y$. $\ZR(Y)$ carries a constant sheaf $K$, a subsheaf of rings $\mathcal{O}$ with stalk $R_v$ at $v$, and a sheaf of totally ordered abelian groups $\Gamma = K^*/\mathcal{O}^*$ such that the induced map of stalks of the quotient map $K^* \to \Gamma$ at $v$ is the valuation $v: K^* \to \Gamma_v$.\footnote{See \cite[Proposition 6.2.3]{SH18}} Let $\Gamma_+ \subset \Gamma$ be the subsheaf of nonnegative sections of $\Gamma$, whose stalk at $v$ is the monoid of nonnegative elements in the value group $\Gamma_v$.

\subsection{Rees algebras} 
We say that $\mathcal{A}$ is a Rees algebra on $Y$ if it is a finitely generated graded $\mathcal{O}_Y$-subalgebra of $\mathcal{O}_Y[T]$.\footnote{An $\mathcal{O}_Y$-algebra $\mathcal{A}$ is a finitely generated graded $\mathcal{O}_Y$-algebra if there is an affine covering $Y=\cup U_i$ such that $\mathcal{A}|_{U_i}$ is the sheaf associated to a finitely generated graded $\mathcal{O}_Y(U_i)$-algebra.} This generalizes the notion of a Rees algebra of an ideal. If $\mathcal{A}=\oplus_{m \geq 0} \mathcal{I}_mT^m$ is a Rees algebra, let $\mathcal{A}_+=\oplus_{m>0}\mathcal{I}_mT^m$. We say that $\mathcal{A}$ is a nonzero Rees algebra if the irrelevant ideal $\mathcal{A}_+$ is nonzero.

\subsection{Valuative $\mathbb{Q}$-ideals, associated Rees algebras, and idealistic exponents}
$\Gamma \otimes_\mathbb{Z} \mathbb{Q}$ is also a sheaf of totally ordered abelian groups.\footnote{If $(A,\leq)$ is a totally ordered abelian group, then $A \otimes_{\mathbb{Z}} \mathbb{Q}$ is also a totally ordered abelian group by the rule $x \otimes \frac{p}{q} \leq x' \otimes \frac{p'}{q'}$ if $q'px\leq qp'x'$ in $A$. Note that the canonical inclusion $A \into A \otimes_{\mathbb{Z}} \mathbb{Q}$ is order-preserving.} Let $\Gamma_{\mathbb{Q}+}$ be the subsheaf of nonnegative sections of $\Gamma \otimes_\mathbb{Z} \mathbb{Q}$. A valuative $\mathbb{Q}$-ideal is a section $\gamma \in H^0(\ZR(Y), \Gamma_{\mathbb{Q}+})$. Every valuative $\mathbb{Q}$-ideal determines an ideal sheaf $\mathcal{I}_{\gamma}$ by 
$
\mathcal{I}_{\gamma}(U)=\{f \in \mathcal{O}_Y(U) | v(f) \geq \gamma_v \forall{v \in \pi_Y^{-1}(U)}\}
$. Furthermore, $\gamma$ determines a graded $\mathcal{O}_Y$-subalgebra of $\mathcal{O}_Y[T]$ by $\mathcal{A}_{\gamma}=\oplus_{m \geq 0} \mathcal{I}_{m \gamma}T^m$. The associated graded algebra $\mathcal{A}_{\gamma}$ is the integral closure of $\mathcal{A}$ in $\mathcal{O}_Y[T,T^{-1}]$.

On the other hand, every nonzero Rees algebra $\mathcal{A}$ determines a valuative $\mathbb{Q}$-ideal $\gamma_{\mathcal{A}}$ by $\gamma_{\mathcal{A},v} = \min\{ \tfrac{1}{n} v(f) | fT^n \in (\mathcal{A}_+)_{\pi_Y(v)}\}$. A valuative $\mathbb{Q}$-ideal $\gamma$ is called an idealistic exponent if $\gamma=\gamma_{\mathcal{A}}$ for some nonzero Rees algebra $\mathcal{A}$. If $\gamma=\gamma_{\mathcal{A}}$ is an idealistic exponent, then the associated graded algebra $\mathcal{A}_{\gamma}$ is a nonzero Rees algebra. This establishes a one-to-one correspondence between idealistic exponents and nonzero Rees algebras integrally closed in $\mathcal{O}_Y[T,T^{-1}]$.

\subsection{Locally representing idealistic exponents} 
If $U \subset Y$ is open and $\gamma$ is a valuative $\mathbb{Q}$-ideal on $Y$, we denote the restriction $\gamma_U := \gamma|_{\pi_Y^{-1}(U)}$. Given an affine open $U \subset Y$, a nonempty finite collection $f_i \in \mathcal{O}_Y(U)$, and $a_i \in \mathbb{Q}_{>0}$, we write
$$
(f_1^{a_1}, \dots, f_r^{a_r}) := ( \min\{ a_i  \cdot v(f_i) \})_v \in H^0(\pi_Y^{-1}(U), \Gamma_{\mathbb{Q}+})
$$
for the naturally associated valuative $\mathbb{Q}$-ideal on $U$.\footnote{Note that $\pi_Y^{-1}(U)=\ZR(U)$.} Note that $(1)$ is the zero valuative $\mathbb{Q}$-ideal. For integers $w_i > 0$, $(f_1^{1/w_1}, \dots, f_{r}^{1/w_r})$ is the idealistic exponent associated to the $\mathcal{O}_U$-subalgebra of $\mathcal{O}_U[T]$ generated by the sections $f_iT^{w_i}$'s. Conversely, let $\mathcal{A}$ be a nonzero Rees algebra on $Y$, with $\gamma=\gamma_{\mathcal{A}}$ its associated idealistic exponent, and suppose $f_1T^{w_1}, \dots, f_rT^{w_r}$ generate $\mathcal{A}(U)$ as an $\mathcal{O}_Y(U)$-algebra. Then $\gamma_U=(f_1^{1/w_1}, \dots, f_r^{1/w_r})$. Thus we can represent idealistic exponents affine-locally in the form $(f_1^{1/w_1}, \dots, f_{r}^{1/w_r})$, and any valuative $\mathbb{Q}$-ideal represented locally on $Y$ as $(f_1^{1/w_1}, \dots, f_r^{1/w_r})$ is an idealistic exponent.

\subsection{Regular parameters, centers, and reduced centers}
\label{subsec: centers}

\begin{definition}[Regular parameters] 
\label{def: regular parameters}
Let $U=\Spec A$ be a smooth irreducible affine variety. We say that a nonempty subset $x_1, \dots, x_r \in A$ are regular parameters on $U$ if the subvariety $V(x_1,\dots,x_r) = \Spec A/(x_1,\dots,x_r)$ is nonempty and $x_1,\dots,x_r$ extend to a regular system of parameters in $\mathcal{O}_{U,p}$ for every $p \in V(x_1,\dots,x_r)$. Equivalently, $V(x_1,\dots,x_r)$ is a nonempty, regular, pure codimension $r$ subvariety of $U$.
\end{definition}

For example, $x,y \in k[x,y,z]$ are regular parameters on $\mathbb{A}^3_k=\Spec k[x,y,z]$.

\begin{definition}[Centers and their associated reduced centers] 
\label{def: centers}
Assume further that $Y$ is a smooth variety. A center on $Y$ is a valuative $\mathbb{Q}$-ideal $J$ on $Y$ for which there is an affine covering $Y= \cup_i U_i$ and $a_1, \dots, a_r \in \mathbb{Q}_{>0}$ such that for each $i$, either $J_{U_i} = (1)$ or $J_{U_i} = (x_{i1}^{a_1}, \dots, x_{ir}^{a_r})$, where $x_{i1}, \dots, x_{ir}$ are regular parameters on $U_i$. Let $w_1,\dots,w_r$ be the unique positive integers such that $(\ell/w_1, \dots, \ell/w_r)=(a_1,\dots,a_r)$ for some $\ell \in \mathbb{Q}$. Define the reduced center $\bar{J}$ associated to $J$ to be the valuative $\mathbb{Q}$-ideal on $Y$ such that for each $i$, $\bar{J}_{U_i} = (x_{i1}^{1/w_1}, \dots, x_{ir}^{1/w_r})$ if $J_{U_i} = (x_{i1}^{a_1}, \dots, x_{ir}^{a_r})$, and $\bar{J}_{U_i} = (1)$ if $J_{U_i}=(1)$.
\end{definition}

For example, $(x^{3/5}, y^{4/5})$ is a center on $\mathbb{A}^3_k=\Spec k[x,y,z]$, and its associated reduced center is $(x^{1/4}, y^{1/3})$.

\subsection{Reducible varieties}
If $Y=\sqcup_i Y_i$ is the disjoint union of finitely many separated integral varieties $Y_i$ over $k$, then the Zariski-Riemann space of $Y$ is $\ZR(Y) := \sqcup_i \ZR(Y_i)$. All notions developed for separated integral varieties naturally extend to this setting. For example, a center on $Y$ is given by a center on each $Y_i$.

\section{Derivatives and Orders of vanishing}
\label{sec: derivatives}

Throughout this section, let $k$ be a field. Unless stated otherwise, no assumptions are put on $k$.

\subsection{Derivatives of ideals}
Let $A$ be a $k$-algebra and $I \subset A$ an ideal. Define  $D^{\leq n}_{A/k} I \subset A$ to be the ideal generated by the image under the evaluation map $\Diff^n_{A/k} \times I \to A$ sending $(D,f) \mapsto D(f)$, where $\Diff^n_{A/k} \subset \Hom_k(A,A)$ is the module of differential operators of order $\leq n$ \cite[16.8.1]{EGAIV}. If $Y$ is a variety over $k$ and $\mathcal{I} \subset \mathcal{O}_Y$ a coherent ideal, then define $\mathcal{D}^{\leq n}_{Y/k}\mathcal{I} \subset \OY$ to be the ideal sheaf generated by the image under the evaluation map $\sDiff^n_{Y/k} \times \mathcal{I} \to \mathcal{O}_Y$, where $\sDiff^n_{Y/k}$ is the $\mathcal{O}_Y$-module of differential operators of order $\leq n$ \cite[16.8.1]{EGAIV}. The ideal $\mathcal{D}^{\leq n}_{Y/k}\mathcal{I}$ is coherent. For an affine open $U \subset Y$, we have $(\mathcal{D}^{\leq n}_{Y/k}\mathcal{I})(U)=D^{\leq n}_{\mathcal{O}(U)/k} \mathcal{I}(U)$, and for a point $p \in Y$, we have $(\mathcal{D}^{\leq n}_{Y/k}\mathcal{I})_p=D^{\leq n}_{\mathcal{O}_{Y,p}/k} \mathcal{I}_p$. If $Y$ is a smooth variety over a perfect field $k$, then $\mathcal{D}^{\leq 1}_{Y/k}\mathcal{I}$ is generated by the image under the evaluation map $\Der_{k}(\mathcal{O}_Y,\mathcal{O}_Y) \times \mathcal{I} \to \mathcal{O}_Y$ sending $(\nabla , f) \mapsto \nabla(f)$, where $\Der_{k}(\mathcal{O}_Y,\mathcal{O}_Y)$ is the sheaf of $k$-derivations from $\mathcal{O}_Y$ to itself. Moreover, if $k$ has characteristic zero, then the higher derivatives can be obtained iteratively by $D^{\leq n}_{Y/k} \mathcal{I} = D^{\leq 1}_{Y/k}(D^{\leq n-1}_{Y/k} \mathcal{I})$.

In the constructive resolution literature, $\Delta(\mathcal{I})$ is common notation for $\mathcal{D}^{\leq 1}_{Y/k}\mathcal{I}$. For example, see \cite[Definition 13.6]{BEV05}. In \cite{Kol07}, Koll{\'a}r writes $D(\mathcal{I})$ for $\mathcal{D}^{\leq 1}_{Y/k}\mathcal{I}$ \cite[Definition 3.73]{Kol07}.

\subsection{Computing derivatives}
After proving the following proposition and corollary, we will provide Algorithm \nameref{Diff}, which explicitly computes derivatives of ideals on smooth varieties. In \cite{BDFPRR18}, the following proposition and corollary are proven in the special case of a smooth complete intersection over an algebraically closed field. 

\begin{proposition}
\label{Kahler basis}
Let $A = k[x_1, \dots, x_N]/(f_1, \dots, f_r)$, and suppose that $\Spec A$ is a smooth variety of pure dimension $n$. Let $J = (\frac{\partial f_i}{\partial x_j})_{ij}$ be the Jacobian; let $ROW \subset \{1, \dots, r\}$ and $\COL \subset \{1, \dots, N\}$ be increasing $(N-n)$-tuples; and let $M =  (\frac{\partial f_i}{\partial x_j})_{i \in ROW, j \in COL}$ be the corresponding submatrix of the Jacobian. Let $h = \det M$,\footnote{If $M$ is the $0 \times 0$ matrix, then its determinant $h=\det M=1$.} and let $C=( C_{ij} )_{i \in ROW,j \in COL}$ be the cofactor matrix of $M$. For each $j' \not\in \COL$, that is for each $j' \in \{1, \dots, N\} \setminus \COL$, define
$$
D_{j'} := h \frac{\partial }{\partial x_{j'}} - \sum_{\substack{i \in \ROW\\ j \in \COL}}   \frac{\partial f_i}{\partial x_{j'}} C_{ij}  \frac{\partial }{\partial x_j}
$$
Then the $D_{j'}$'s form a basis for $\Der_k(A[h^{-1}], A[h^{-1}])$.
\end{proposition}
\proof
Denote $df := \sum_{\ell=1}^N \frac{\partial f}{\partial x_{\ell}} dx_{\ell}$, so that $J = [df_1 \cdots df_r]^T$. We first observe that the submodule $(df_i)_{i \in \ROW}$ in $\oplus_{\ell=1}^N A[h^{-1}] dx_{\ell}$ is free with basis $\{ hdx_j + \sum_{j' \not\in \COL, i \in \ROW} C_{ij} \frac{\partial f_i}{\partial x_{j'}}  dx_{j'}\}_{j \in \COL}$ because $C^T M = hI$, where $I$ is the identity matrix. Thus we have isomorphisms
\begin{align*}
\Omega_{A[h^{-1}]/k}
&= \bigoplus\limits_{\ell=1}^N A[h^{-1}] dx_{\ell} \bigg/ (df_1,\dots,df_r) \\
&= \bigoplus\limits_{\ell=1}^N A[h^{-1}] dx_{\ell} \bigg/ \bigg( (df_{i'})_{i' \not\in \ROW} + (hdx_j + \textstyle\sum_{j' \not\in \COL, i \in \ROW} C_{ij} \frac{\partial f_i}{\partial x_{j'}}  dx_{j'})_{j \in \COL} \bigg)   \\
&= \bigoplus\limits_{j' \not\in \COL} A[h^{-1}] dx_{j'} \bigg/ (\widetilde{df_{i'}})_{i' \not\in \ROW}
\end{align*}
where the last isomorphism substitutes $dx_j = -h^{-1}\sum_{j' \not\in \COL, i \in \ROW} C_{ij} \frac{\partial f_i}{\partial x_{j'}}  dx_{j'}$ for every $j \in \COL$. The stalk $(\Omega_{A[h^{-1}]/k})_{\frak{p}}$ is free of rank $n$ for every $\frak{p} \in \Spec A[h^{-1}]$, which means that the coefficients of the $\widetilde{df_{i'}}$'s must be contained in every prime ideal of $A[h^{-1}]$. But $A[h^{-1}]$ is smooth over $k$, in particular reduced, hence $(\widetilde{df_{i'}})_{i' \not\in \ROW}=0$, and $\Omega_{A[h^{-1}]/k} = \oplus_{j' \not\in \COL} A[h^{-1}] dx_{j'}$. We are done once we see that the derivation in $\Der_k(A[h^{-1}], A[h^{-1}])=\Hom_{A[h^{-1}]}(\Omega_{A[h^{-1}]/k}, A[h^{-1}])$ corresponding to the coordinate projection on $dx_{j'}$ is $h^{-1}D_{j'}$. To see this, just observe the following computation:
\begin{align*}
h df &=  h \sum_{j' \not\in \COL} \frac{\partial f}{\partial x_{j'}} dx_{j'} + h \sum_{j \in \COL} \frac{\partial f}{\partial x_j} dx_j \\
&= h \sum_{j' \not\in \COL} \frac{\partial f}{\partial x_{j'}} dx_{j'} - \sum_{j \in \COL} \sum_{j' \not\in \COL} \frac{\partial f}{\partial x_j} \bigg( \sum_{i \in \ROW} C_{ij} \frac{\partial f_i}{\partial x_{j'}}  \bigg)dx_{j'} \\
&= \sum_{j' \not\in \COL} \bigg( h \frac{\partial f}{\partial x_{j'}} - \sum_{\substack{i \in \ROW\\j \in \COL}} \frac{\partial f_i}{\partial x_{j'}} C_{ij} \frac{\partial f}{\partial x_j} \bigg)dx_{j'}      \\
&= \sum_{j' \not\in \COL}  D_{j'}( f) dx_{j'}
\end{align*} \qed

\begin{corollary} 
\label{diff generators}
Let the notation be as in Proposition \ref{Kahler basis}, and let $I \subset A$ be an ideal generated by the images of polynomials $g_1,\dots,g_m \in k[x_1,\dots,x_N]$. Then the extension of the following ideal
$$
(g_1,\dots,g_m)+\big( D_{j'}( g_1), \dots,  D_{j'}( g_m) \big)_{j' \not\in \COL} \subset k[x_1,\dots,x_N]
$$
in $A[h^{-1}]$ is equal to $D^{\leq 1}_{A[h^{-1}]/k}(I)$. 
\end{corollary}
\proof By Proposition \ref{Kahler basis}, $\Der_k(A[h^{-1}],A[h^{-1}])$ is a free $A[h^{-1}]$-module with basis $\{ D_{j'} \}_{j' \not\in \COL}$. Thus $\Diff^1_{A[h^{-1}]/k}$ is a free $A[h^{-1}]$-module with basis $\{id_{A[h^{-1}]} \} \cup \{ D_{j'} \}_{j' \not\in \COL}$, where $id_{A[h^{-1}]}$ is the identity on $A[h^{-1}]$. Because $D_{j'}( \sum_i h_i g_i) = \sum_i D_{j'}(h_i)id_{A[h^{-1}]}(g_i) + h_i D_{j'}(g_i)$, the image under the evaluation map $\Diff^1_{A[h^{-1}]/k} \times I \to A[h^{-1}]$ is generated by $\{id_{A[h^{-1}]}(g_i), D_{j'}( g_i) \}_{i=1,j' \not\in \COL}^m$. \qed

\subsection{Algorithm for taking derivatives of ideals}
The author learned about the algorithm for taking derivatives of ideals from \cite[Algorithm 2 \texttt{DeltaCheck}, page 16]{BDFPRR18} and \cite[\texttt{Delta}]{FP19}. Algorithm \nameref{Diff} is the same as \cite[Algorithm Delta, page 312]{F07}, but we add further explanations on both the theoretical and computational aspects involved in computing derivatives of ideals. 
\Suppressnumber
\begin{lstlisting}[caption={Diff},label={Diff}]
$\text{\bf{Input:}}$ $f_1,\dots,f_r, g_1, \dots, g_m \in k[x_1,\dots,x_N]$, such that $\Spec A$ is a smooth variety over $k$ of pure dimension where $A = k[x_1,\dots,x_N]/(f_1,\dots,f_r)$. Let $I \subset A$ be the ideal generated by the images of $g_1,\dots,g_m$.
$\text{\bf{Output:}}$ An ideal in $k[x_1,\dots,x_N]$ whose extension in $A$ equals $D^{\leq 1}_{A/k}I$. 
$\text{\bf{Remark:}}$ If $k$ has characteristic zero, then this algorithm can iteratively compute higher derivatives by $D^{\leq r}_{A/k}I=D^{\leq 1}_{A/k}(D^{\leq r-1}_{A/k}I)$. Also, this algorithm can extend to computing derivatives of ideals on smooth varieties not of pure dimension: by using Algorithm @\nameref{orthogonal idempotents}@, we can compute the derivative of an ideal on each component of a smooth variety, then glue them together using Proposition @\ref{gluing ideals}@. @\Reactivatenumber@ 
Initialization
  @{\bu}@$n=\dim A$ (so $\Spec A$ is of pure dimension $n$)@\footnote{Use {\sc Singular}'s procedure \texttt{dim} \cite[5.1.23]{DGPS19} to compute Krull dimension.}@
  @{\bu}@$J=[df_1 \cdots df_r]^T = (\frac{\partial f_i}{\partial x_j} )_{ij}$
  @{\bu}@@$L = \{ \text{$N-n$ by $N-n$ square submatrices of the Jacobian $J$} \}$@
  @{\bu}@$I_{return} = (1) \subset k[x_1,\dots,x_N]$, so $I_{return}=k[x_1,\dots,x_N]$
Because $\Spec A$ is smooth over $k$ of pure dimension $n$, we have @$$\Spec A = \bigcup_{M \in L} \Spec A[\det M^{-1}]$$@
We may replace $L$ by a subset for which $\Spec A = \bigcup_{M \in L} \Spec A[\det M^{-1}]$ still holds.@\footnote{See {\sc Singular}'s procedure \texttt{lift} \cite[5.1.75]{DGPS19}.}@ 
$\text{\bf{for}}$ each $M \in L$ $\text{\bf{do}}$
  As in the notation of Proposition @\ref{Kahler basis}@, let 
    @{\bu}@$h=\det M$; 
    @{\bu}@$\ROW \subset \{1,\dots,r\}$ be the row indices and $\COL \subset \{1,\dots,N\}$ be the column indices of the Jacobian $J$ that its submatrix $M$ involve; 
    @{\bu}@$C$ be the cofactor matrix of $M$;
    @{\bu}@For every $j' \not\in \COL$, $D_{j'}=h \frac{\partial }{\partial x_{j'}} - \sum_{\substack{i \in \ROW\\ j \in \COL}}   \frac{\partial f_i}{\partial x_{j'}} C_{ij}  \frac{\partial }{\partial x_j}$.
  Set @$$I_M=(f_1,\dots,f_r)+(g_1,\dots,g_m)+\big(D_{j'}(g_1), \dots, D_{j'}(g_m)\big)_{j' \not\in \COL} \subset k[x_1,\dots,x_N]$$@
  Then let $I_{return}=I_{return} \cap (I_M : h^{\infty})$
By Corollary @\ref{diff generators}@, Proposition @\ref{gluing ideals}@, and computational remark @\ref{subsec: compute remark}@, we get that the extension of the ideal $I_{return}$ in $A$ equals $D^{\leq 1}_{A/k}I$. 
$\text{\bf{return}}$ $I_{return}$
\end{lstlisting}

\subsection{Example for Diff}
For polynomials $F_1,\dots,F_N \in k[x_1,\dots,x_n]$, we have 
$$
D^{\leq 1}_{k[x_1,\dots,x_n]/k} (F_1,\dots,F_N) = ( \tfrac{ \partial F_i}{\partial x_j} )_{ij}
$$ 
by Algorithm \nameref{Diff}.

\subsection{Order of vanishing}
Let $Y$ be a smooth variety over a field $k$ of characteristic zero, and let $\mathcal{I} \subset \OY$ be a coherent ideal on $Y$. Define the order of vanishing of $\mathcal{I}$ at a point $p \in Y$ to be
$$
\ord_p \mathcal{I} := \max\{ r \geq 0 \ | \  \mathcal{I}_p \subset \frak{m}_p^r \}
$$
Define the maximal order of vanishing of $\mathcal{I}$ on $Y$ to be
$$
\maxord_Y \mathcal{I} = \maxord \mathcal{I} := \max \{ \ord_p \mathcal{I} \ | \ p \in Y \}
$$
The order of vanishing of $\mathcal{I}$ is a function $p \mapsto \ord_p \mathcal{I}$ on points $p \in Y$ taking values in $\mathbb{N} \cup \{\infty\}$. Note that $\ord_p \mathcal{I} \geq 1$ iff $p \in V(\mathcal{I})$, and that $\ord_p \mathcal{I} = \infty$ iff $\mathcal{I}_p=0$ iff $\mathcal{I}$ vanishes on the irreducible component of $Y$ containing $p$. Thus $\maxord \mathcal{I} < \infty$ iff $\mathcal{I}$ does not vanish on any irreducible component of $Y$. By going to the completions of the stalks of $Y$, we can see that the locus of points where $\mathcal{I}$ vanishes to order at least $b$ is cut out by the ideal $\mathcal{D}^{\leq b-1}_{Y/k} \mathcal{I}$, leading us towards the following definition.

\begin{definition}[$b$-singular locus] Let $b \in \mathbb{N}$. Define the $b$-singular locus $V(\mathcal{I},b)$ to be the closed subvariety of $Y$ cut out by $\mathcal{D}^{\leq b-1}_{Y/k} \mathcal{I}$.
\end{definition}

The $b$-singular locus $V(\mathcal{I},b)$ is precisely the locus of points where $\mathcal{I}$ vanishes to order at least $b$, that is, $p \in V(\mathcal{I},b)$ iff $\ord_p \mathcal{I} \geq b$. Note that we can explicitly compute the $b$-singular locus by iteratively applying Algorithm \nameref{Diff} $b-1$ times to the ideal $\mathcal{I}$.

\subsection{Algorithm for maximal order of vanishing} 
We provide the algorithm for computing the maximal order of vanishing of any ideal on a smooth variety below. We note that the procedure \texttt{DeltaList} in the {\sc Singular} library \texttt{resolve.lib} computes the maximal order of vanishing of ideals with finite maximal order of vanishing on smooth varieties of pure dimension \cite[\texttt{DeltaList}]{FP19}.
\Suppressnumber
\begin{lstlisting}[caption={maximalOrderOfVanishing},label={maximalOrderOfVanishing}]
Let $k$ be a field of characteristic zero.
$\text{\bf{Input:}}$ $f_1,\dots,f_r, g_1, \dots, g_m \in k[x_1,\dots,x_N]$, such that $\Spec A$ is a smooth variety over $k$, where $A = k[x_1,\dots,x_N]/(f_1,\dots,f_r)$. Let $I \subset A$ be the ideal generated by the images of $g_1,\dots,g_m$.
$\text{\bf{Output:}}$ The maximal order of vanishing $\maxord I$ of the coherent ideal $I$ on $\Spec A$. 
$\text{\bf{Remark:}}$ Let $Y$ be a smooth variety over $k$ and $\mathcal{I} \subset \OY$ a coherent ideal on $Y$. Suppose we have a covering of $Y$ by finitely many open affines $Y= \cup_i U_i$. By this algorithm, we can compute $\maxord_{U_i} \mathcal{I}|_{U_i}$ for each $i$. Then the maximal order of vanishing of $\mathcal{I}$ on $Y$ is $\maxord_Y \mathcal{I} = \max_i \{ \maxord_{U_i} \mathcal{I}|_{U_i}\}$. @\Reactivatenumber@ 
Let $\frak{p}_i \subset k[x_1,\dots,x_N]$ be the associated primes of the ideal $(f_1,\dots,f_r) \subset k[x_1,\dots,x_N]$.@\footnote{See {\sc Singular}'s primary decomposition library \cite{PDSL19}.}@ 
Because $\Spec A$ is a smooth variety, the associated primes $\{\frak{p}_i\}_i$ are the minimal primes of $A$ and we have the irreducible component decomposition $\Spec A = \sqcup_i \Spec A/\frak{p}_i$. 
$\text{\bf{if}}$ $(g_1,\dots,g_m) \subset \frak{p}_i$ for some $i$ $\text{\bf{then}}$
  $I$ vanishes on the irreducible component of $\Spec A$ corresponding to $\frak{p}_i$.
  $\text{\bf{return}}$ $\maxord I = \infty$
$\text{\bf{else}}$
  In this case, $I$ doesn't vanish on any irreducible component of $\Spec A$, so $\maxord I < \infty$.
  By repeatedly applying Algorithm @\nameref{Diff}@, obtain the smallest integer $b \geq 0$ such that $D^{\leq b}_{A/k}I=(1)$. Let $b$ be this integer. 
  $\text{\bf{return}}$ $\maxord I = b$
\end{lstlisting}

\subsection{Example for maximalOrderOfVanishing}
Consider the coherent ideal $\mathcal{I}=(z^2-x^2y^2)$ on $\mathbb{A}^3_k=\Spec k[x,y,z]$. By repeatedly applying Algorithm \nameref{Diff}, we get the following:
$$
D^{\leq 1}_{k[x,y,z]/k}(z^2-x^2y^2)  = (-2xy^2, -2x^2y, 2z) = (xy^2, x^2y, z)
$$
$$
D^{\leq 2}_{k[x,y,z]/k}(z^2-x^2y^2) = D^{\leq 1}_{k[x,y,z]/k}(xy^2, x^2y, z) = (1)
$$
Thus $\maxord \mathcal{I} = 2$. 


\section{Maximal contact hypersurfaces}
\label{sec: maximal contact}
Let $Y$ be a smooth variety over a field $k$ of characteristic zero, and let $\mathcal{I} \subset \mathcal{O}_Y$ be a coherent ideal with $b=\maxord \mathcal{I}<\infty$. A smooth hypersurface\footnote{See Appendix \ref{subsec: smooth hypersurface}} $H \subset Y$ is a maximal contact hypersurface if $H$ scheme-theoretically contains $V(\mathcal{I},b)$, that is, $\mathcal{D}^{\leq b-1}_{Y/k}\mathcal{I}$ contains the ideal cutting $H$ out of $Y$. If $U \subset Y$ is an open affine, then a smooth hypersurface $H \subset U$ is a local maximal contact hypersurface of $\mathcal{I}$ on $U$ if $H$ scheme-theoretically contains $V(\mathcal{I},b)|_U$. 

Because $k$ has characteristic zero, there exists a cover $Y = \cup_i U_i$ by finitely many affine opens $U_i$ such that there is a local maximal contact hypersurface $H_i \subset U_i$ of $\mathcal{I}$ on $U_i$ for each $i$. Note that if $U_i$ does not meet $V(\mathcal{I},b)$, then the empty subvariety of $U_i$, which is a smooth hypersurface of $U_i$, contains $V(\mathcal{I},b)|_{U_i} = \emptyset$.

\subsection{Algorithm for maximal contact}
We provide below the algorithm for obtaining hypersurfaces of maximal contact. We note that the {\sc Singular} library \texttt{resolve.lib} has an incorrect algorithm for obtaining maximal contact hypersurfaces. For most examples of practical interest, however, \texttt{resolve.lib} does correctly compute hypersurfaces of maximal contact, and this is mostly because the maximal contact algorithm in \texttt{resolve.lib} and Algorithm \nameref{maximalContact} both agree up to line 5 of Algorithm \nameref{maximalContact}.\footnote{See the static procedure \texttt{Coeff} in the {\sc Singular} library \texttt{resolve.lib} \cite{FP19}. Static procedures can only be viewed in the source code of {\sc Singular} libraries.} The problem with \texttt{resolve.lib} lies in its static procedure \texttt{CoverCenter}, which corresponds to what comes after line 5 in Algorithm \nameref{maximalContact}. Though \texttt{CoverCenter} works correctly in the special case of computing maximal contact hypersurfaces on affine space, \texttt{CoverCenter} incorrectly computes maximal contact hypersurfaces on general smooth varieties. Nevertheless, the ideas found in \texttt{CoverCenter} are still valid, and Algorithm \nameref{maximalContact} uses these ideas to provide the correct algorithm for computing maximal contact hypersurfaces.

\Suppressnumber
\begin{lstlisting}[caption={maximalContact},label={maximalContact}]
Let $k$ be a field of characteristic zero.
$\text{\bf{Input:}}$ $f_1,\dots,f_r \in k[x_1,\dots,x_N]$ and $I \subset A$ an ideal, where  $A = k[x_1,\dots,x_N]/(f_1,\dots,f_r)$, such that $Y=\Spec A$ is a smooth variety of pure dimension over $k$. Let $\mathcal{I} \subset \mathcal{O}_Y$ be the coherent ideal associated to $I$. Assume that $0 < \maxord \mathcal{I} < \infty$. 
$\text{\bf{Output:}}$ A finite list $\MC=\{ (G_i, F_i) \}_i$, where $F_i,G_i \in k[x_1,\dots,x_N]$, such that
@{\bu}@$\Spec A=\bigcup_i \Spec A[G_i^{-1}]$, where each $\Spec A[G_i^{-1}]$ is nonempty
@{\bu}@For each $i$, $F_i$ is a local maximal contact hypersurface of $\mathcal{I}$ on $\Spec A[G_i^{-1}]$
@{\bu}@For each $i$, $F_i=1$ iff $V(\mathcal{I},b) \cap \Spec A[G_i^{-1}] = \emptyset$ 
@{\bu}@$\MC$ is a list of size one iff there exists $(G,F) \in \MC$ such that $G=1$. @\Reactivatenumber@ 
Let $b=\maxord \mathcal{I}$.
Use Algorithm @\nameref{Diff}@ to obtain $h_1,\dots,h_m \in k[x_1,\dots,x_N]$ whose images generate $D^{\leq b-1}_{A/k}I$.
Use Algorithm @\nameref{isSmoothHypersurface}@ to determine whether any $h_i$ cuts out a smooth hypersurface out of $Y$, since in that case $h_i$ would be a global maximal contact hypersurface of $\mathcal{I}$ on $Y$.
$\text{\bf{if}}$ there is an $i$ such that $\Spec A/h_i$ is a smooth hypersurface of $\Spec A$ $\text{\bf{then}}$
  $\text{\bf{return}}$ $\MC=\{ (1,h_i) \}$
Initialization
  @{\bu}@$\MC=\emptyset$
  @{\bu}@$n=\dim A$
  @{\bu}@$J=[df_1 \cdots df_r]^T = (\frac{\partial f_i}{\partial x_j} )_{ij}$
  @{\bu}@@$L = \{ \text{$N-n$ by $N-n$ square submatrices of the Jacobian $J$} \}$@
We may replace $L$ by a subset for which $\Spec A=\bigcup_{M \in L} \Spec A[\det M^{-1}]$ still holds.
Use Algorithm @\nameref{orthogonal idempotents}@ to obtain $e_1, \dots, e_d \in k[x_1, \dots, x_N]$, the orthogonal idempotents of $A$. This means $\Spec A$ has $d$ irreducible components, and its $t$th component is $\Spec A[e_t^{-1}]$.
$\text{\bf{for}}$ each $M \in L$ $\text{\bf{do}}$
  As in the notation of Proposition @\ref{Kahler basis}@, let 
    @{\bu}@$h=\det M$; 
    @{\bu}@$\ROW \subset \{1,\dots,r\}$ be the row indices and $\COL \subset \{1,\dots,N\}$ be the column indices of the Jacobian $J$ that its submatrix $M$ involve; 
    @{\bu}@$C$ be the cofactor matrix of $M$;
    @{\bu}@For every $j' \not\in \COL$, $D_{j'}=h \frac{\partial }{\partial x_{j'}} - \sum_{\substack{i \in \ROW\\ j \in \COL}}   \frac{\partial f_i}{\partial x_{j'}} C_{ij}  \frac{\partial }{\partial x_j}$.
  Note that $\Spec A[h^{-1}]=\bigsqcup_{t=1}^d \Spec A[h^{-1},e_t^{-1}]$, and that the irreducible components of $\Spec A[h^{-1}]$ correspond to the $e_t$ for which $\Spec A[h^{-1},e_t^{-1}]$ is nonempty.
  Introducing a new variable $y$, we express $A[h^{-1},e_t^{-1}]=A[(he_t)^{-1}]$ as follows: @$$A[h^{-1},e_t^{-1}] = \frac{k[x_1, \dots,x_N, y]}{(f_1, \dots,f_r, 1- yhe_t)}$$@
  Note that $\Spec A[h^{-1},e_t^{-1}]$ is nonempty iff $1 \not\in (f_1, \dots,f_r, 1- yhe_t)$.@\footnote{See {\sc Singular}'s procedure \texttt{reduce} \cite[5.1.124]{DGPS19}}@
  $\text{\bf{for}}$ each $e_t \in \{e_1,\dots,e_d\}$ such that $\Spec A[h^{-1},e_t^{-1}]$ is nonempty $\text{\bf{do}}$
    Because $b=\maxord \mathcal{I}$, by Corollary @\ref{diff generators}@, we have that @$$(h_1,\dots,h_m) + \big( D_{j'}(h_1), \dots, D_{j'}(h_m) \big)_{j' \not\in \COL}+(f_1, \dots,f_r, 1- yhe_t)=(1) \subset k[x_1,\dots,x_N,y]$$@
    Thus we can find polynomials $a_i,b_{ij'} \in k[x_1,\dots,x_N,y]$@\footnote{Use {\sc Singular}'s procedure \texttt{lift} \cite[5.1.75]{DGPS19} to compute the $a_i,b_{ij'}$'s}@ such that @$$\sum_{i=1}^m  a_i h_i +  \sum_{j' \not\in \text{COL}} \sum_{i=1}^m b_{ij'} D_{j'}(h_i) \equiv 1 \mod (f_1, \dots,f_r, 1- yhe_t)$$@
    This means that we have the covering @$$\Spec A[h^{-1},e_t^{-1}] = \bigcup_{\substack{a_i \neq 0 \text{ in} \\ A[h^{-1},e_t^{-1}]}} \Spec A[h^{-1},e_t^{-1},h_i^{-1}] \ \cup  \bigcup_{\substack{b_{ij'} \neq 0 \text{ in} \\A[h^{-1},e_t^{-1}]}}  \Spec A[h^{-1},e_t^{-1},D_{j'}(h_i)^{-1}] $$@
    If $h_i=0$ in $A[h^{-1},e_t^{-1}]$, set $a_i=0 \in k[x_1,\dots,x_N,y]$, and if $D_{j'}(h_i)=0$ in $A[h^{-1},e_t^{-1}]$, set $b_{ij'}=0 \in k[x_1,\dots,x_N,y]$. Note that the equation in line 23 and the equality in line 24 still hold. We make these changes so that the following statement is true: if $a_i \neq 0$ in $A[h^{-1},e_t^{-1}]$, then $\Spec A[h^{-1},e_t^{-1},h_i^{-1}]$ is nonempty, and if $b_{ij'} \neq 0$ in $A[h^{-1},e_t^{-1}]$, then $\Spec A[h^{-1},e_t^{-1},D_{j'}(h_i)^{-1}]$ is nonempty. 
  $\text{\bf{for}}$ each $i$ such that $a_i \neq 0$ in $A[h^{-1},e_t^{-1}]$ $\text{\bf{do}}$
    $\Spec A[h^{-1}, e_t^{-1},h_i^{-1}]$ does not meet $V(\mathcal{I},b)$. 
    $\MC=\MC,\{ (h\cdot e_t \cdot h_i, 1) \}$
  $\text{\bf{for}}$ each $i$ and $j'$ such that $b_{ij'} \neq 0$ in $A[h^{-1}]$ $\text{\bf{do}}$
    By @\cite[\href{https://stacks.math.columbia.edu/tag/07PF}{Tag 07PF}]{Stacks20}@, $h_i$ cuts a smooth subvariety out of $\Spec A[h^{-1}, e_t^{-1},D_{j'}(h_i)^{-1}]$.
    Because $\Spec A[h^{-1}, e_t^{-1},D_{j'}(h_i)^{-1}]$ is nonempty, $D_{j'}(h_i) \neq 0$ in $A[h^{-1}, e_t^{-1},D_{j'}(h_i)^{-1}]$, hence $h_i \neq 0$ in $A[h^{-1}, e_t^{-1},D_{j'}(h_i)^{-1}]$. Since $A[h^{-1}, e_t^{-1},D_{j'}(h_i)^{-1}]$ is an integral domain, by Krull's principal ideal theorem, $h_i$ cuts out a pure codimension one subvariety out of $\Spec A[h^{-1}, e_t^{-1},D_{j'}(h_i)^{-1}]$.
    Therefore, $h_i$ cuts out a smooth hypersurface out of $\Spec A[h^{-1}, e_t^{-1},D_{j'}(h_i)^{-1}]$. 
    $\text{\bf{if}}$ $(h_1,\dots,h_m) \neq (1)$ in $A[h^{-1}, e_t^{-1},D_{j'}(h_i)^{-1}]$ $\text{\bf{then}}$
      $\Spec A[h^{-1}, e_t^{-1},D_{j'}(h_i)^{-1}]$ meets the $V(\mathcal{I},b)$.
      $\MC=\MC,\{ (h\cdot e_t \cdot D_{j'}(h_i), h_i) \}$
    $\text{\bf{else}}$
      $\Spec A[h^{-1}, e_t^{-1},D_{j'}(h_i)^{-1}]$ does not meet the $V(\mathcal{I},b)$.
      $\MC=\MC,\{ (h\cdot e_t \cdot D_{j'}(h_i), 1) \}$
$\text{\bf{if}}$ there exists $(G,F) \in \MC$ such that $G$ is invertible in $A$ $\text{\bf{then}}$
   $\MC=\{(1,F)\}$
$\text{\bf{return}}$ $\MC$
\end{lstlisting}

\subsection{Example for maximalContact}
Consider the coherent ideal $\mathcal{I}=(x^5+x^3y^3+y^8)$ on the affine plane $\mathbb{A}^2_k=\Spec k[x,y]$. By repeatedly applying Algorithm \nameref{Diff}, we find that $D^{\leq 4}_{k[x,y]/k}\mathcal{I} = (x,y^2)$ and $D^{\leq 5}_{k[x,y]/k}\mathcal{I}=(1)$, so that $\maxord \mathcal{I} = 5$. Since $x \in D^{\leq 4}_{k[x,y]/k}\mathcal{I}$ cuts a smooth hypersurface out of $\mathbb{A}^2_k$, we see that $x$ is the equation of a global maximal contact hypersurface of $\mathcal{I}$ on $\mathbb{A}^2_k$. Thus Algorithm \nameref{maximalContact}, when applied to the ideal $\mathcal{I}$, would return the list $\MC=\{ (1,x) \}$.

\subsection{Algorithm for lifting maximal contact} 
In order to compute (and even to define\footnote{See Definition \ref{lex order}.}) the centers in weighted resolution, it is necessary to lift maximal contact hypersurfaces of an ideal on a smooth subvariety to the ambient variety. We provide below the algorithm for lifting maximal contact hypersurfaces.

\Suppressnumber
\begin{lstlisting}[caption={liftMaximalContact},label={liftMaximalContact}]
Let $k$ be a field of characteristic zero.
$\text{\bf{Input:}}$ $f_1,\dots,f_r,h_1,\dots,h_m \in k[x_1,\dots,x_N]$ and $I \subset A/(h_1,\dots,h_m)$ an ideal, where $A = k[x_1,\dots,x_N]/(f_1,\dots,f_r)$, such that $Y=\Spec A$ and $Z= \Spec A/(h_1,\dots,h_m)$ are smooth varieties of pure dimension over $k$. Let $\mathcal{I} \subset \mathcal{O}_Z$ be the coherent ideal on $Z$ associated to $I$. Assume that $0 < \maxord_Z \mathcal{I} < \infty$. 
$\text{\bf{Output:}}$ A finite list $\MC_{\text{lift}}=\{ (G_i, F_i) \}_i$, where $F_i,G_i \in k[x_1,\dots,x_N]$, such that
@{\bu}@$\Spec A=\bigcup_i \Spec A[G_i^{-1}]$, where each $\Spec A[G_i^{-1}]$ is nonempty
@{\bu}@For each $i$, if $(\Spec A[G_i^{-1}]) \big|_Z$ is nonempty, then $F_i$ is a local maximal contact hypersurface of $\mathcal{I}$ on $(\Spec A[G_i^{-1}]) \big|_Z$
@{\bu}@For each $i$, $F_i=1$ iff $V(\mathcal{I},b) \cap (\Spec A[G_i^{-1}]) \big|_Z = \emptyset$ 
@{\bu}@$\MC_{\text{lift}}$ is a list of size one iff there exists $(G,F) \in \MC_{\text{lift}}$ such that $G=1$. @\Reactivatenumber@ 
Initialization
  @{\bu}@$\MC_{\text{lift}}=\emptyset$
Let $\MC$ be the output of Algorithm @\nameref{maximalContact}@ applied to the polynomials $f_1, \dots, f_r, h_1, \dots, h_m \in k[x_1,\dots,x_N]$ and ideal $I \subset k[x_1, \dots, x_N]/(f_1, \dots, f_r, h_1, \dots, h_m)$. $\MC$ provides the data of local maximal contact hypersurfaces of $\mathcal{I}$ on $Z$.
$\text{\bf{for}}$ each $(G_i,F_i) \in \MC$ $\text{\bf{do}}$
  $\MC_{\text{lift}}= \MC_{\text{lift}}, \{ (G_i,F_i)\}$
The complement of $Z$ in $Y$ is covered by $\Spec A[h_1^{-1}], \dots, \Spec A[h_m^{-1}]$, so the restrictions of these distinguished opens to $Z$ are empty and hence do not meet the $V(\mathcal{I},b)$.
$\text{\bf{for}}$ each $i$ such that $\Spec A[h_i^{-1}]$ is nonempty $\text{\bf{do}}$
  $\MC_{\text{lift}}= \MC_{\text{lift}}, \{ (h_i,1)\}$
$\text{\bf{if}}$ there exists $(G,F) \in \MC_{\text{lift}}$ such that $G$ is invertible in $A/(h_1,\dots,h_m)$ $\text{\bf{then}}$
   $\MC_{\text{lift}}=\{(1,F)\}$
$\text{\bf{return}}$ $\MC_{\text{lift}}$
\end{lstlisting}

\subsection{Example for liftMaximalContact}
Let $Y=\Spec k[x,y,z]$ and $Z =V(x,y) \subset Y$. Consider the coherent ideal $\mathcal{I}=(z^5)$ on $Z$. Applying Algorithm \nameref{maximalContact} to $\mathcal{I}$ returns the list $\MC=\{ (1,z) \}$. Thus applying Algorithm \nameref{liftMaximalContact} to $\mathcal{I}$, with ambient variety $Y$, returns the list $\MC_{\text{lift}}=\{ (1,z) \}$.


\section{Lexicographic order invariant and the associated center}
\label{sec: lexinv}

Throughout this section, let $Y$ be a smooth variety over a field $k$ of characteristic zero, and $\mathcal{I} \subset \OY$ a coherent ideal on $Y$.

\subsection{Coefficient ideals}
Let $b \in \mathbb{N}$. Define the coefficient ideal
$$
\mathcal{C}_e(\mathcal{I},b) := \sum_{i=0}^{b-1} \big( \mathcal{D}^{\leq i}_{Y/k} \mathcal{I} \big)^{b!/(b-i)}
$$
In \cite{ATW20}, Abramovich, Temkin, and W{\l}odarczyk use the following larger coefficient ideal\footnote{See \cite[Definition 4.1.1]{ATW20}.}
$$
\mathcal{C}(\mathcal{I},b) := \sum_{\sum_{i=1}^{b-1} (b-i) c_i \geq b!    } (\mathcal{D}^{\leq 0}_{Y/k} \mathcal{I})^{c_0}  \cdot  (\mathcal{D}^{\leq 1}_{Y/k} \mathcal{I})^{c_1}   \cdots    (\mathcal{D}^{\leq b-1}_{Y/k}\mathcal{I})^{c_{b-1}} 
$$
Our coefficient ideal $\mathcal{C}_e(\mathcal{I},b)$ is contained in $\mathcal{C}(\mathcal{I},b)$, and most times this containment is proper. However, there is no issue with replacing $\mathcal{C}(\mathcal{I},b)$ by $\mathcal{C}_e(\mathcal{I},b)$ in Abramovich-Temkin-W{\l}odarczyk's paper \cite{ATW20}: in particular, the invariant defined in \cite[Section 5.1]{ATW20} is unchanged when computing with our smaller coefficient ideal. The larger coefficient ideal $\mathcal{C}(\mathcal{I},b)$ is useful for theoretical proofs, but it is computationally more efficient to handle the smaller coefficient ideal $\mathcal{C}_e(\mathcal{I},b)$ (the subscript $e$ stands for ``efficient").

\subsection{The ordered sets $\mathbb{Q}^{\leq n}_{\geq 0}$ and $\mathbb{N}_{>0}^{\leq n}$} 
Let $\mathbb{Q}_{\geq 0}$ denote the nonnegative rational numbers, and $\mathbb{N}_{>0}$ the positive natural numbers. Define
$
\mathbb{Q}_{\geq 0}^{\leq n} := \sqcup_{i=0}^n \mathbb{Q}_{\geq 0}^i
$
as the disjoint union of the products $\mathbb{Q}_{\geq 0}^i$ for $i=0,\dots,n$. Order elements of $\mathbb{Q}_{\geq 0}^{\leq n}$ lexicographically, with truncated sequences considered larger. For example, 
$$
(1, 1, 1) < (1, 1, 2) < (1, 2, 1) < (1, 2) < (2, 2, 1) < ()
$$
In this way, $\mathbb{Q}_{\geq 0}^{\leq n}$ becomes totally ordered. Also, define
$
\mathbb{N}_{>0}^{\leq n} := \sqcup_{i=0}^n \mathbb{N}_{>0}^i
$
with lexicographic order, just like $\mathbb{Q}_{\geq 0}^{\leq n}$, where truncated sequences are considered larger. In this way, $\mathbb{N}_{>0}^{\leq n}$ becomes a well-ordered set.

\subsection{Lexicographic order invariant, $b$-invariant, and associated parameters}
Let the dimension of $Y$ be $n$. We will define a function called the $b$-invariant with respect to $\mathcal{I}$ taking values in a well-ordered set:
\begin{align*}
\binv \mathcal{I}: V(\mathcal{I}) &\to \mathbb{N}_{>0}^{\leq n}  \\
p &\mapsto  \binv_p \mathcal{I}
\end{align*}
From this we will define the lexicographic order invariant:
\begin{align*}
\lexinv \mathcal{I}: V(\mathcal{I}) &\to \mathbb{Q}_{\geq 0}^{\leq n}  \\
p &\mapsto  \lexinv_p \mathcal{I} 
\end{align*}

\begin{definition}[$\btoa$ function] 
\label{b to a}
The function $\btoa: \mathbb{N}_{>0}^{\leq n} \to \mathbb{Q}_{\geq 0}^{\leq n}$ is defined recursively by sending the empty sequence $()$ to itself and sending $(b_1,...,b_r) \in \mathbb{N}_{>0}^{\leq n}$ to $(a_1,\dots, a_r) \in \mathbb{Q}_{\geq 0}^{\leq n}$, where $a_1=b_1$ and $(a_2,\dots,a_r) = \btoa(b_2,\dots,b_r)/(b_1-1)!$. 
\end{definition}

For example, $\btoa( 3, 8,  50400) = (3,4,5)$. Note that the function $\btoa$ is order preserving and injective, hence the image $\btoa (\mathbb{N}_{>0}^{\leq n}) \subset \mathbb{Q}_{\geq 0}^{\leq n}$ is a well-ordered subset.

\begin{definition}[$b$-invariant, Lexicographic order invariant, and associated parameters, {see \cite{ATW20}}] 
\label{lex order}
Let $p \in V(\mathcal{I})$, and set $\mathcal{I}[1] :=\mathcal{I}$ and $b_1 := \ord_{\mathcal{O}_{X,p}}(\mathcal{I}[1])$. If $b_1 = \infty$, define $\binv_p \mathcal{I} = ()$ to be the empty sequence and associate the empty sequence of parameters. Otherwise, choose an affine neighborhood $U_1 \subset U$ of $p$ with a function $x_1 \in \mathcal{O}(U_1)$ such that $x_1$ is the equation of a local maximal contact hypersurface of $\mathcal{I}[1]$ on $U_1$. Let $V(x_1) := \Spec \mathcal{O}(U_1)/x_1$. Inductively, write
\begin{align*}
\mathcal{I}[i+1] &:= \mathcal{C}_e(\mathcal{I}[i] , b_i)\big|_{V(x_1,\dots,x_i)}
\end{align*}
Note that $\mathcal{I}[i+1]$ is a coherent ideal on $V(x_1,\dots,x_i)$. Set 
$$
b_{i+1} := \ord_p \mathcal{I}[i+1]
$$ 
If $b_{i+1}=\infty$, define $\binv_p \mathcal{I} = (b_1, \dots, b_i)$ and associate the sequence of parameters $(x_1,\dots,x_i)$ to the point $p$. Otherwise, choose an affine neighborhood $U_{i+1} \subset U_i$ of $p$ with a function $x_{i+1} \in \mathcal{O}(U_{i+1})$ such that $x_{i+1}$ is the equation of a local maximal contact hypersurface of $\mathcal{I}[i+1]$ on $U_{i+1}|_{V(x_1,\dots,x_i)}$.\footnote{Really, $x_{i+1}$ is the lift of the equation of a local maximal contact hypersurface of $\mathcal{I}[i+1]$ on $U_{i+1}|_{V(x_1,\dots,x_i)}$.} Let $V(x_1,\dots,x_{i+1}) = \Spec \mathcal{O}(U_{i+1})/(x_1,\dots,x_{i+1})$.

Suppose that $\binv_p\mathcal{I}=(b_1,\dots,b_r)$ and that $x_1,\dots,x_r \in \mathcal{O}(U_r)$ are the sequence of parameters associated to $p$.\footnote{Note that $x_1,\dots,x_r$ are regular parameters on $U_r$.} Define the lexicographic order invariant 
$$
\lexinv_p \mathcal{I}  := (\btoa)( \binv_p \mathcal{I}  )=\btoa( b_1,\dots,b_r)
$$
and call the sequence $(x_1,\dots,x_r)$ the associated parameters of $\mathcal{I}$ at $p$. If $X=V(\mathcal{I})$, we also write $\lexinv_p X:=\lexinv_p \mathcal{I}$, keeping in mind the embedding $X \subset Y$. 
\end{definition}

The lexicographic order invariant $\inv_p$ is independent of choices and upper-semicontinuous by \cite[Theorem 5.1.1]{ATW20}. On the other hand, the associated parameters are not well-defined and evidently depend on choices.

\begin{definition}[Maximal lexicographic order invariant] 
Define the maximal lexicographic order invariant of $\mathcal{I}$ on $Y$ to be
$$
\maxinv_Y \mathcal{I} = \maxinv \mathcal{I} := \max \{  \inv_p \mathcal{I} \ | \ p \in V(\mathcal{I}) \}
$$
If $X=V(\mathcal{I})$, then we also write $\maxinv X:=\maxinv \mathcal{I}$, keeping in mind the embedding $X \subset Y$. 
\end{definition}

The invariant $\inv_p$ detects singularities. The subvariety $V(\mathcal{I})$ is regular at $p \in V(\mathcal{I})$ iff $\inv_p \mathcal{I}$ is equal to a constant sequence of ones. If $V(\mathcal{I})$ is pure codimension $c$ in $Y$, then $V(\mathcal{I})$ is nonsingular iff $\maxinv \mathcal{I} = (1,\dots,1)$, the constant sequence of length $c$. The locus of points attaining the maximal lexicographic order invariant is a closed subset of $Y$: these are the points on $V(\mathcal{I})$ with the ``worst" singularities, as measured by the invariant $\inv_p$. 

The parameters associated to a point attaining the maximal lexicographic order invariant locally cut out the locus of maximal lexicographic order invariant. Due to the importance of this fact, we state this in the following proposition.

\begin{proposition} Let the notation be as in Definition \ref{lex order}. Then we have that
$$
\lexinv_p \mathcal{I}  = \min_{q \in V(x_1,\dots,x_r) } \lexinv_q \mathcal{I} 
$$
In particular, if $\lexinv_p \mathcal{I} = \maxinv \mathcal{I}$, then $V(x_1,\dots,x_r)$ is precisely the locus of points inside $U_r$ attaining the maximal lexicographic order invariant of $\mathcal{I}$.
\end{proposition}
\proof Follows easily by induction. \qed

\subsection{The associated center}
Suppose that $\maxinv \mathcal{I} = (a_1,\dots,a_r)$. Let $Y = \cup_i U_i$ be a cover by affine opens such that following holds. If $U_i$ contains a point attaining the maximal lexicographic order invariant, let $x_{i1},\dots,x_{ir} \in \mathcal{O}(U_i)$ be such that $(x_{i1},\dots,x_{ir})$ is the sequence of associated parameters of $\mathcal{I}$ at every point of $U_i$ in the vanishing locus of $x_{i1},\dots,x_{ir}$. Let $J$ be the center on $Y$ such that if $U_i$ contains a point attaining the maximal lexicographic order invariant, then $J_{U_i} = (x_{i1}^{a_1}, \dots, x_{ir}^{a_r})$, and otherwise $J_{U_i} = (1)$. By \cite[Theorem 5.3.1]{ATW20}, $J$ is indeed a valuative $\mathbb{Q}$-ideal, since the local idealistic exponents must glue by uniqueness. We call $J$ the unique center associated to $\mathcal{I}$. The reduced center $\bar{J}$ associated to $J$ is the center associated to the pair $(V(\mathcal{I}) \subset Y)$ in \cite[Theorem 1.1.1]{ATW20}. 

In the following section, we provide the algorithm that explicitly computes the maximal lexicographic order invariant and the associated parameters, which is precisely the data of the center associated to $\mathcal{I}$.

\subsection{Algorithm for the associated center} 
We now show how to explicitly compute the maximal lexicographic order invariant and the associated parameters on some open cover. The algorithm proceeds as a competition among open sets covering the smooth variety. Each descent in dimension to a local hypersurface of maximal contact corresponds to each round of the competition. The remaining contenders after a round are those open sets that attain a truncation of the maximal lexicographic order invariant with length the codimension corresponding to the round. The contenders that are dropped after a round are designated as losers. The competition finishes once the maximal lexicographic invariant is found, in which case the surviving contenders are winners, while all other past contenders are losers. 
\Suppressnumber
\begin{lstlisting}[caption={prepareCenter},label={prepareCenter}]
Let $k$ be a field of characteristic zero.
$\text{\bf{Input:}}$ Smooth affine varieties $Y_1,\dots,Y_N$ over $k$, and coherent ideals $\mathcal{I}_i \subset \mathcal{O}_{Y_i}$.@\footnote{The data of the affine variety $Y_1$, for example, is inputted as a polynomial ring $R$ with an ideal $I\subset R$ such that $\mathcal{O}(Y_1)=R/I$. The data of $\mathcal{I}_1 \subset \mathcal{O}_{Y_1}$ is given by an ideal in $R$ whose extension in $R/I$ corresponds to $\mathcal{I}_1$.}@ Let $Y = \sqcup_i Y_i$ and $\mathcal{I} \subset \OY$ be the coherent ideal such that $\mathcal{I}|_{Y_i}= \mathcal{I}_i$. 
$\text{\bf{Compute:}}$ The maximal lexicographic order invariant and associated parameters of $\mathcal{I}$.
$\text{\bf{Output:}}$ $\maxinv \mathcal{I}=(a_1,\dots,a_r)$ and finite lists $\win=\{ (U_i, [x_{i1}, \dots, x_{ir}], \ell_i ) \}_i$ and $\lose=\{ (W_j, \ell_j) \}_j$ such that
@{\bu}@Each $U_i$ is a nonempty distinguished affine open of $Y_{\ell_i}$, and similarly, each $W_j$ is a nonempty distinguished affine open of $Y_{\ell_j}$. 
@{\bu}@For each $i$, there exists a point $p \in U_i$ such that $\inv_p \mathcal{I} = \maxinv \mathcal{I}$.
@{\bu}@For each $j$, $\inv_p \mathcal{I} < \maxinv \mathcal{I}$ at every point $p \in W_j$.
@{\bu}@$x_{i1}, \dots, x_{ir} \in \mathcal{O}(U_i)$, and $(x_{i1}, \dots, x_{ir})$ is the sequence of associated parameters of $\mathcal{I}$ at each point $p \in U_i$ such that $\inv_p \mathcal{I} = \maxinv \mathcal{I}$.
$\text{\bf{Remark:}}$ The center $J$ associated to $\mathcal{I}$ is given by $J|_{U_i} = (x_{i1}^{a_1}, \dots, x_{ir}^{a_r})$ for each $i$, and $J|_{W_j}=0$ for each $j$. We obtain the reduced center $\bar{J}$ associated to $J$ as follows. Suppose $a_i=n_i/m_i$ with $n_i,m_i \in \mathbb{N}_{>0}$. Set $w_i = \lcm(da_1, \dots,da_r)/ da_i$, where $d=m_1 \cdots m_r$. Then $\bar{J}_{U_i} = (x_{i1}^{1/w_1}, \dots, x_{ir}^{1/w_r})$ for each $i$, and $\bar{J}_{W_j}=(1)$ for each $j$.@\footnote{See \cite[Theorem 1.1.1]{ATW20}}@ @\Reactivatenumber@ 
Initialization
  @{\bu}@$\maxbinv= ()$
  @{\bu}@$\text{maxinvfound}=\F$
  @{\bu}@$\con=\{ (Y_1, \mathcal{I}_1, [ \ ], 1), \dots, (Y_N, \mathcal{I}_N, [ \ ], N) \}$
  @{\bu}@$\win=\emptyset$
  @{\bu}@$\lose=\emptyset$
$\text{\bf{while}}$ $\text{NOT maxinvfound}$ $\text{\bf{do}}$
  $\bmax = -\infty$
  $\text{\bf{for}}$ each $(U, \mathcal{I}_Z ,[h_1,\dots,h_m], \ell) \in \con$ $\text{\bf{do}}$
    @{\bu}@$h_1,\dots,h_m \in \mathcal{O}(U)$ are regular parameters on the affine open $U \subset Y$.
    @{\bu}@Let $Z := V(h_1,\dots,h_m)=\Spec \mathcal{O}(U)/(h_1,\dots,h_m)$.
    @{\bu}@$\mathcal{I}_Z \subset \mathcal{O}_Z$ is a coherent ideal on $Z$.
    $\bcurr=\maxord \mathcal{I}_Z$, the maximal order of vanishing of $\mathcal{I}_Z$ on $Z$.
    $\bmax=\max( \bmax, \bcurr)$
  $\surv=\emptyset$
  $\text{\bf{for}}$ each $(U, \mathcal{I}_Z ,[h_1,\dots,h_m], \ell) \in \con$ $\text{\bf{do}}$
    $\text{\bf{if}}$ $\maxord \mathcal{I}_Z=\bmax$ $\text{\bf{then}}$
      $\surv=\surv,\{ (U, \mathcal{I}_Z ,[h_1,\dots,h_m], \ell) \}$
    $\text{\bf{else}}$
      $\lose=\lose,\{ (U, \ell) \}$
  $\text{\bf{if}}$ $\bmax=\infty$ $\text{\bf{then}}$
    The maximal lexicographic order invariant has been found! 
    All survivors are now winners. 
    $\text{maxinvfound}=\F$
    $\win=\win,\surv$
  $\text{\bf{else}}$
    $\maxbinv=\maxbinv,\bmax$
    $\con=\emptyset$
      $\text{\bf{for}}$ each $(U, \mathcal{I}_Z ,[h_1,\dots,h_m], \ell) \in \surv$ $\text{\bf{do}}$
        @{\bu}@$Z := V(h_1,\dots,h_m)$
        @{\bu}@$\mathcal{I}_Z \subset \mathcal{O}_Z$ is a coherent ideal on $Z$.
        Apply Algorithm @\ref{liftMaximalContact}@ to the closed subvariety $Z \subset U$, obtaining finite lists $\{ (U_i, f_i) \}_i$ and $\{W_j\}_j$ such that 
          @{\bu}@each $U_i$ and $W_j$ are nonempty distinguished affine opens of $U$
          @{\bu}@$U = (\cup_i U_i) \cup (\cup_j W_j)$
          @{\bu}@For each $i$, $V(\mathcal{I}_Z, \bmax) \cap U_i|_Z \neq \emptyset$
          @{\bu}@For each $j$, $V(\mathcal{I}_Z, \bmax) \cap W_j|_Z = \varnothing$
          @{\bu}@$f_i$ is the equation of a local maximal contact hypersurface of $\mathcal{I}_Z$ on $U_i|_Z$.
        For each $i$, let $V(f_i) := \Spec \mathcal{O}(U_i)/(h_1,\dots,h_m,f_i)$.
        $\con=\con,\{ ( U_i, \mathcal{C}_e( \mathcal{I}_Z, \bmax)|_{V(f_i)}, [h_1,\dots,h_m,f_i], \ell) \}_i$
        $\lose=\lose, \{(W_j, \ell) \}_j$
Note that $\maxinv\mathcal{I} = (\btoa)(\maxbinv)$.
$\text{\bf{return}}$ $(\btoa)(\maxbinv), \win, \lose$
\end{lstlisting}

\subsection{Example for prepareCenter}
We use Algorithm \nameref{prepareCenter} to compute the maximal lexicographic order invariant and associated parameters of the ideal $\mathcal{I}[1]=\mathcal{I}=(x^5+x^3y^3+y^8)$ on $Y=\Spec k[x,y]$. Since $\mathcal{D}^{\leq 4}_{Y/k}\mathcal{I}[1] = (x,y^2)$, we have that $b_1=\maxord \mathcal{I}[1]=5$, and we take $x$ as the first hypersurface of maximal contact. A direct computation finds that 
$$
\mathcal{I}[2]=\mathcal{C}_e(\mathcal{I}[1],5)|_{V(x)}=(y^{180})
$$
Thus $b_2=180$, and $y$ is the second hypersurface of maximal contact. Thus the maximal $b$-invariant of $\mathcal{I}$ is $(5,180)$, and the sequence of associated parameters is $(x,y)$. This means that $\maxinv\mathcal{I} = \btoa(5,180)=(5,15/2)$, hence the center $J$ associated to $\mathcal{I}$ is $(x^5, y^{15/2})$, with reduced center $\bar{J}=(x^{1/3},y^{1/2})$.


\section{Weighted blowing ups}
\label{sec: blowup}

We begin with the definition of weighted blowing ups.

\begin{definition}[Weighted blowing ups, see {\cite[\S 3.3]{ATW20} and \cite[\S 4]{Que20}}] 
Let $\gamma$ be a center on a smooth variety $Y$, and let $\mathcal{A}_{\gamma} \subset \mathcal{O}_Y[T]$ be its associated Rees algebra. We define the (stack-theoretic) weighted blowings up of $Y$ along the center $\gamma$ to be 
$$
Bl_Y(\gamma) :=\cProj_Y \mathcal{A}_{\gamma} = [(\Spec_{\OY} \mathcal{A}_{\gamma}  \setminus S_0) /\mathbb{G}_m]
$$
where the vertex $S_0$ is the closed subscheme defined by the irrelevant ideal $\mathcal{A}_{\gamma,+}=\oplus_{m>0}\mathcal{I}_{m\gamma}T^m$ and the grading on $\mathcal{A}_{\gamma}$ defines a $\mathbb{G}_m$-action $(T,s) \mapsto T^m \cdot s$ for $s \in \mathcal{I}_{m\gamma}$. 
\end{definition}

Since weighted blowing ups can be computed affine-locally, we may assume that the smooth variety $Y$ is irreducible and affine, and that the center $\gamma$ on $Y$ is either $(1)$ or $(x_1^{a_1}, \dots,x_r^{a_r})$, where $x_1,\dots,x_r\in A$ are regular parameters on $Y$ and $a_1,\dots,a_r \in \mathbb{Q}_{>0}$. If $\gamma=(1)$, then $Bl_Y(\gamma)=\cProj_Y \mathcal{O}_Y[T]=Y$, so we pay attention to the case when $\gamma=(x_1^{a_1}, \dots,x_r^{a_r})$. We assume that the weights $a_i$ are the inverse of positive integers $w_i$, since the weighted resolution algorithm only blows up along centers locally of the form $(x_1^{1/w_1}, \dots,x_r^{1/w_r})$.

\subsection{The setup}
For the rest of this section, let $Y=\Spec A$ be a smooth irreducible affine variety over a field $k$ of characteristic zero, $x_1,\dots, x_r \in A$ regular parameters on $Y$,\footnote{See Definition \ref{def: regular parameters}.} and $w_1,\dots,w_r$ positive integers. By the end of this section we will provide an algorithm for computing the weighted blowings up of $Y$ along the center $\gamma =(x_1^{1/w_1},\dots, x_r^{1/w_r})$.

Let $\mathcal{A}_{\gamma}=\oplus_{m \geq 0} \mathcal{I}_{m \gamma}T^m \subset \mathcal{O}_Y[T]$ be the Rees algebra associated to $\gamma$, and let $I_{m \gamma} = \mathcal{I}_{m \gamma}(Y)$ and $A_{\gamma}=\mathcal{A}_{\gamma}(Y)=\oplus_{m \geq 0}I_{m \gamma}T^m \subset A[T]$ be the global sections. Let 
$$
A^{pre}_{\gamma} := A[x_1T^{w_1}, \dots, x_rT^{w_r}] \subset A[T]
$$
The valuative $\mathbb{Q}$-ideal on $Y$ associated to the Rees algebra $A^{pre}_{\gamma}$ is $\gamma$, hence $A_{\gamma}$ is the integral closure of $A^{pre}_{\gamma}$ in $A[T,T^{-1}]$.

\subsection{The generators of the Rees algebra associated to a center}

\begin{proposition}[See {\cite[Section 3.4]{ATW20}}]
\label{Rees algebra generators}
For every $m\geq 0$, 
$$
I_{m \gamma} = \big( x_1^{b_1} \cdots x_r^{b_r}  \ \big| \  \sum w_ib_i \geq m \big)
$$
In particular, $A_{\gamma}=\oplus_{m\geq 0} I_{m\gamma}T^m \subset A[T]$ is generated by $\{ x_iT,x_iT^2, \dots,x_iT^{w_i} \}_{i=1}^r$ as an $A$-algebra.
\end{proposition}
\proof
Let $b_1,\dots,b_r$ be nonnegative integers such that $\sum w_ib_i \geq m$. Then for every $v \in \ZR(Y)$, we have
$$
v(x_1^{b_1} \cdots x_r^{b_r})=\sum b_i v(f_i) = \sum b_i w_i \tfrac{1}{w_i}v(x_i) \geq \sum w_ib_i \gamma_v \geq m\gamma_v
$$
thus obtaining the inclusion $( x_1^{b_1} \cdots x_r^{b_r}  \ | \  \sum w_ib_i \geq m ) \subset I_{m \gamma}$.

For the reverse inclusion, note that the nonempty subvariety $V(x_1,\dots,x_r)=\Spec A/(x_1,\dots,x_r)$ is regular and of pure dimension $n-r$. Let $y \in Y$ be the generic point of some irreducible component of $V(x_1,\dots,x_r)$. Then $x_1,\dots,x_r$ form a regular system of parameters in $\mathcal{O}_{Y,y}$. Let $\frak{m}_y$ be the maximal ideal of $\mathcal{O}_{Y,y}$ and $\kappa(y)=\mathcal{O}_{Y,y}/\frak{m}_y$ be the residue field, and let $\hat{\mathcal{O}}_{Y,y}$ be the $\frak{m}_y$-adic completion of $\mathcal{O}_{Y,y}$. By the Cohen structure theorem, there is a ring map $\kappa(y) \to \hat{\mathcal{O}}_{Y,y}$ (not necessarily a $k$-algebra map) such that $\kappa(y) \to \hat{\mathcal{O}}_{Y,y} \to \hat{\mathcal{O}}_{Y,y}/\hat{\frak{m}}_y$ is an isomorphism; by this ring map we fix a $\kappa(y)$-algebra structure on $\hat{\mathcal{O}}_{Y,y}$. Again by the Cohen structure theorem, $\hat{\mathcal{O}}_{Y,y}=\kappa(y)[[ x_1, \dots,x_r]]$ is a formal power series ring in $r$ variables $x_1,\dots, x_r$ because $\mathcal{O}_{Y,y}$ is a regular local ring, where $x_i \in \hat{\mathcal{O}}_{Y,y}$ is the image of $x_i \in \mathcal{O}_{Y,y}$. 

Now consider the monomial valuation $v: \kappa(y)[[ x_1, \dots,x_r]] \to \mathbb{Z} \cup \{ \infty\}$ given on variables by $v(x_i)=w_i$ and extended to all power series by
$$
v\bigg(\sum_{s \in \mathbb{N}^r} a_s x^s \bigg) = \min \bigg\{ \sum s_i v(x_i) | a_s \neq 0 \bigg\}
$$
where $s_i$ is the $i$th component of the multi-index $s \in \mathbb{N}^r$, $a_s \in \kappa(y)$, and $x^s:=x_1^{s_1} \cdots x_r^{s_r}$ (see \cite[Definition 6.1.4]{SH18}). Restrict $v$ to a valuation on $\mathcal{O}_{Y,y}$, which then extends uniquely to a valuation on the function field $K$ of $Y$. The valuation $v$ vanishes on $k^*$ because the invertible elements in $\kappa(y)[[ x_1, \dots,x_r]]$ are precisely the formal power series with nonzero constant term. Observe that $\gamma_v = 1$ and $v$ is nonnegative on $\mathcal{O}_{Y,y}$. Because $v$ is strictly positive on the regular parameters $x_1, \dots, x_r$, we get that $v$ has center $y$ on $Y$, hence $v \in \ZR(Y)$.

Let $f \in A$ not be in the ideal $(x_1^{b_1} \cdots x_r^{b_r}  \ \big| \  \sum w_ib_i \geq m )$. Because $\mathcal{O}_{Y,y} \into \hat{\mathcal{O}}_{Y,y}$ is faithfully flat, the image of $f$ in $\hat{\mathcal{O}}_{Y,y}$ is not contained in the ideal generated by the monomials $\prod x_i^{b_i}$ such that $\sum w_ib_i \geq m$ (see \cite[Theorem 7.5]{Mat89}). Thus the image of $f$ in $\kappa(y)[[ x_1, \dots,x_r]]$ must be a power series involving some monomial $\prod x_i^{b_i}$ with $\sum w_ib_i < m$. We then obtain $v(f) < m = m \gamma_v$, hence $f \not\in I_{m \gamma}$ and we are done showing that $I_{m \gamma} = \big( x_1^{b_1} \cdots x_r^{b_r}  \ \big| \  \sum w_ib_i \geq m \big)$. 

To see that $\{ x_iT,x_iT^2, \dots,x_iT^{w_i} \}_{i=1}^r$ generates $A_{\gamma}$ as an $A$-algebra, first observe that $I_{m\gamma} = \big( x_1^{b_1} \cdots x_r^{b_r}  \ \big| \  \sum b_i \leq m \leq \sum w_ib_i \big)$. After all, if $m < \sum b_i$, then there are nonnegative integers $a_i\leq b_i$ such that $m=\sum (b_i-a_i)$, hence 
$$
\prod x_i^{b_i}=\prod x_i^{a_i} \prod x_i^{b_i-a_i} \in \big( x_1^{b_1} \cdots x_r^{b_r}  \ \big| \  \sum b_i \leq m \leq \sum w_ib_i \big)
$$
Now, suppose $\sum b_i \leq m \leq \sum w_ib_i$. Then there are integers $b_i \leq m_i \leq w_ib_i$ such that $\sum m_i=m$, and then there are integers $1\leq w_{i,1}, \dots, w_{i,b_i} \leq w_i$ such that $w_{i,1}+\cdots+w_{i,b_i}=m_i$. Thus we have
$$
x_1^{b_1} \cdots x_r^{b_r}T^m = \prod_{i=1}^r x_i^{b_i} T^{m_i} =\prod_{i=1}^r (x_iT^{w_{i,1}} \cdots x_iT^{w_{i,b_i}})
$$
and we get that $A_{\gamma}$ is generated by $\{ x_iT,x_iT^2, \dots,x_iT^{w_i} \}_{i=1}^r$ as an $A$-algebra. \qed

\subsection{The charts of the weighted blowings up} The weighted blowings up of $\gamma$ on $Y$ is 
$$
Bl_{Y}(\gamma) := [ (\Spec A_{\gamma} \setminus S_0) / \mathbb{G}_m ]
$$
Let $y_i := x_iT^{w_i} \in A_{\gamma}$ for each $i$. By Proposition \ref{Rees algebra generators}, we can see that the radical of the ideal in $A_{\gamma}$ generated by the $y_i$'s contains the irrelevant ideal $S_0$, so that $\Spec A_{\gamma} \setminus S_0$ is covered by the distinguished opens $\Spec A_{\gamma}[y_i^{-1}]$. Since $y_i \in A_{\gamma}$ is homogeneous, the $\mathbb{G}_m$-action on $\Spec A_{\gamma}$ restricts to each $\Spec A_{\gamma}[y_i^{-1}]$. Thus $Bl_{Y}(\gamma)$ is covered by the open substacks $[\Spec A_{\gamma}[y_i^{-1}]/\mathbb{G}_m]$. We define $[\Spec A_{\gamma}[y_i^{-1}]/\mathbb{G}_m]$ to be the $x_i$-chart of $Bl_{Y}(\gamma)$. Let $W_i:=\Spec A_{\gamma}[y_i^{-1}]/(y_i-1)$. Since $y_i-1$ has degree zero when restricting the $\mathbb{Z}$-grading on $A_{\gamma}[y_i^{-1}]$ to a $\mathbb{Z}/w_i$-grading, $W_i$ is stablized by the subgroup $\bmu_{w_1} \subset \mathbb{G}_m$. By \cite[Section 3.5]{ATW20} (or \cite[Lemma 4.4.1]{Que20}), the closed embedding $W_i \into \Spec A_{\gamma}[y_i^{-1}]$ gives rise to an isomorphism $[W_i/\bmu_{w_i}] \to [\Spec A_{\gamma}[y_i^{-1}]/\mathbb{G}_m]$.\footnote{Forming the quotient stack $[W_i/\bmu_{w_i}]$ relies on the characteristic of $k$ being zero; if the characteristic was positive and $w_i$ divided the characteristic, then $\bmu_{w_i}$ would not be a smooth group scheme.} Thus $W_i$ is an {\'e}tale cover of the $x_i$-chart, and we also call $W_i$ the $x_i$-chart. Note that the induced morphism
$$
\bigsqcup_{i=1}^r W_i \to Bl_{Y}(\gamma)
$$
is surjective and {\'e}tale.

The map $A \to A_{\gamma}[y_i^{-1}]/(y_i-1)$ that gives $A_{\gamma}[y_i^{-1}]/(y_i-1)$ its $A$-algebra structure defines the composition $W_i \to [W_i/\bmu_{w_i}] \to Y$. Note that the image of $A \to A_{\gamma}[y_i^{-1}]/(y_i-1)$ lies in the zero graded piece with respect to the $\mathbb{Z}/w_i$-grading on $A_{\gamma}[y_i^{-1}]/(y_i-1)$, hence indeed the morphism $W_i \to Y$ descends to $[W_i/\bmu_{w_i}] \to Y$. 

In the following subsections we will show that $\mathcal{O}(W_i) = A_{\gamma}[y_i^{-1}]/(y_i-1)$ is an integral smooth $k$-algebra with the same dimension as $A$, as well as describe the algorithm for its explicit computation.

\subsection{Useful representations of the charts}
There are several useful ways to describe the coordinate ring of the charts $W_i$. Without loss of generality, we focus on the $x_1$-chart $W_1$. First note the canonical isomorphism 
$$
A_{\gamma}[y_1^{-1}]/(y_1-1) = A_{\gamma}/(y_1-1)
$$
We will also want a way to express the coordinate ring of $W_1$ using $A^{pre}_{\gamma}$.

\begin{proposition} 
\label{charts are integral}
$x_1T^{w_1}-1$ is a prime element of $A[T]$, and $A_{\gamma}/(x_1T^{w_1}-1)$ and $A^{pre}_{\gamma}/(x_1T^{w_1}-1)$ are subrings of the integral domain $A[T]/(x_1T^{w_1}-1)$. In particular, the $x_1$-chart $W_1$ is integral.
\end{proposition}
\proof 
Let $\frak{p} \in \Spec A/(x_1,\dots,x_r)$. We first show that $x_1T^{w_1}-1$ is a prime element in $A_{\frak{p}}[T]$. Since $A_{\frak{p}}$ is a regular local ring, in particular a UFD, $A_{\frak{p}}[T]$ is a UFD, hence it suffices to show that $x_1T^{w_1}-1$ is an irreducible element in $A_{\frak{p}}[T]$. Suppose $x_1T^{w_1}-1=(\sum_{i=1}^s a_iT^i)(\sum_{j=1}^{\ell} b_jT^j)$, where $\sum a_iT^i, \sum b_jT^j \in A_{\frak{p}}[T]$. Then $x_1=a_s b_{\ell}$, hence because $x_1$ is irreducible in $A_{\frak{p}}$, we may assume without loss of generality and up to multiplying by a unit of $A_{\frak{p}}$ that $a_s=x_1$ and $b_{\ell}=1$. If $s=w_1$, then $\sum a_iT^i = x_1T^{w_1}-1$. So assume that $s < w_1$, so that $\ell>0$. Then by matching coefficients
\begin{gather*}
a_{s-1}+x_1 b_{\ell-1}=0 \\
a_{s-2}+a_{s-1}b_{\ell-1}+x_1b_{\ell-2}=0 \\
\cdots \\
a_0+a_1b_{\ell-1}+\cdots=0
\end{gather*}
which implies that $a_{s-1},\dots, a_1,a_0 \in x_1A_{\frak{p}}$. But $a_0b_0=-1$ contradicts that $a_0 \in x_1A_{\frak{p}}$. Thus $x_1T^{w_1}-1$ is a prime element in $A_{\frak{p}}[T]$. 

Now let $R \subset A_{\frak{p}}[T]$ be a graded $A$-subalgebra containing $x_1T^{w_1}-1$. We claim that
$$
R \cap \big( (x_1T^{w_1}-1)A_{\frak{p}}[T] \big) = (x_1T^{w_1}-1)R
$$
To obtain this, we must show that if $F \in A_{\frak{p}}[T]$ is a polynomial such that $(x_1T^{w_1}-1)F \in R$, then $F \in R$. First define the length of a polynomial $\sum a_i T^i$ to be $\max \{i | a_i \neq 0\} - \min \{ i | a_i \neq 0\}$. If the length of $F$ is less than $w_1$, then $F$ is equal to a truncation of $-(x_1T^{w_1}-1)F$, hence $F \in R$ because $R$ is a graded subring of $A[T]$. Otherwise, if the length of $F$ is greater than or equal to $w_1$, then the lowest degree term, say $a_mT^m$, of $F$ is also the lowest degree term in $(1-x_1T^{w_1})F$, hence $a_mT^m \in R$, and $(x_1T^{w_1}-1)(F-a_mT^m) \in R$ since $x_1T^{w_1}-1 \in R$. The length of $F-a_mT^m$ is less than the length of $F$, so by inducting on length, we can conclude that $F \in R$. 

Thus $R/(x_1T^{w_1}-1)$ is a subring of $A_{\frak{p}}[T]/(x_1T^{w_1}-1)$, and $R/(x_1T^{w_1}-1)$ is an integral domain. Setting $R=A[T]$, we get that $A[T]/(x_1T^{w_1}-1)$ is an integral domain. By setting $R = A_{\gamma}$ and using that $A_{\gamma} \subset A[T]$, we obtain
$$
A_{\gamma} \cap \big( (x_1T^{w_1}-1)A[T] \big) =A_{\gamma} \cap A[T] \cap \big( (x_1T^{w_1}-1)A_{\frak{p}}[T] \big) = (x_1T^{w_1}-1)A_{\gamma}
$$
hence $A_{\gamma}/(x_1T^{w_1}-1)$ is a subring of the integral domain $A[T]/(x_1T^{w_1}-1)$. The case is exactly the same for $A^{pre}_{\gamma}$. \qed 

\begin{corollary} 
\label{compute coordinate ring}
The $A$-subalgebras $A^{pre}_{\gamma}[T^{-1}]/(x_1T^{w_1}-1)$ and $A_{\gamma}/(x_1T^{w_1}-1)$ of the integral domain $A[T,T^{-1}]/(x_1T^{w_1}-1)$ are equal. In particular, we have the isomorphism $W_1 = \Spec A^{pre}_{\gamma}[T^{-1}]/(x_1T^{w_1}-1)$.
\end{corollary}
\proof By Proposition \ref{charts are integral}, the $A$-algebras $A^{pre}_{\gamma}[T^{-1}]/(x_1T^{w_1}-1)$ and $A_{\gamma}/(x_1T^{w_1}-1)$ are indeed $A$-subalgebras of $A[T,T^{-1}]/(x_1T^{w_1}-1)$. To see equality, it suffices to show that these $A$-subalgebras share the same set of generators. First note that $T^{-1}=x_1T^{w_1-1}$ in $A[T,T^{-1}]/(x_1T^{w_1}-1)$. Thus $T^{-1} \in A_{\gamma}/(x_1T^{w_1}-1)$. By Proposition \ref{Rees algebra generators}, $\{ x_iT,x_iT^2, \dots,x_iT^{w_i} \}_{i=1}^r$ generate $A_{\gamma}$ as an $A$-algebra, so that $T^{-1}, x_2T^{w_2}, \dots, x_3T^{w_3}$ generate $A_{\gamma}/(x_1T^{w_1}-1)$ as an $A$-algebra. But these are precisely the generators of $A^{pre}_{\gamma}[T^{-1}]/(x_1T^{w_1}-1)$, hence we are done. \qed

\subsection{The charts are integral and smooth}\

\begin{lemma} 
\label{charts are affine}
If $A=k[x_1,\dots,x_n]$ is a polynomial ring in $n$ variables and $\gamma=(x_1^{1/w_r},\dots, x_r^{1/w_r})$ is the center on $Y=\mathbb{A}^n_k$, then $A_{\gamma}/(y_1-1)$ is also a polynomial ring in the following $n$ variables 
$$
u, y_2,\dots,y_r,x_{r+1},\dots,x_n
$$
where $u:=T^{-1}=x_1T^{w_1-1}$.
\end{lemma}
\proof
By Proposition \ref{Rees algebra generators}, $A_{\gamma}$ is the $k[x_1,\dots,x_n]$-subalgebra of $k[x_1,\dots,x_n][T]$ generated by $\{x_iT,\dots,x_iT^{w_i}\}_{i=1}^r$. Thus $A_{\gamma}$ is the $k$-subalgebra of $k[x_1,\dots,x_n][T]$ generated by $\{x_i,x_iT,\dots,x_iT^{w_i}\}_{i=1}^r$. $A_{\gamma}/(y_1-1)$ is generated by $u, y_2,\dots,y_r,x_{r+1},\dots,x_n$ as a $k$-algebra because $T^{-1} \in A_{\gamma}/(y_1-1)$. To see that $u, y_2,\dots,y_r,x_{r+1},\dots,x_n$ are algebraically independent over $k$, consider the distinguished open $D(x_1) := \Spec A_{\gamma}[x_1^{-1}]/(y_1-1) \into W_1$. Observe that 
$$
A_{\gamma}[x_1^{-1}]/(y_1-1)= k[x_1,\dots,x_n,T][x_1^{-1}]/(x_1T^{w_1}-1)
$$
Thus $D(x_1)$ has dimension $n$ by Krull's principal ideal theorem. Since $W_1$ is integral by Lemma \ref{charts are integral}, the transcendence degree of the fraction field of $A_{\gamma}/(y_1-1)$ over $k$ must be $n$, hence $u, y_2,\dots,y_r,x_{r+1},\dots,x_n \in A_{\gamma}/(y_1-1)$ must be algebraically independent over $k$, as desired. 

\begin{proposition} The chart $W_1$ is an irreducible smooth variety with the same dimension as $Y$. 
\end{proposition}
\proof 
By Proposition \ref{charts are integral}, $W_1$ is integral, so it suffices to show that $W_1$ is regular. 

Let $p \in \Spec A/(x_1,\dots,x_r)$ be a closed point and $n = \dim Y$. Let $x_{r+1},\dots,x_n \in \mathcal{O}_{Y,p}$ be such that $x_1,\dots,x_n$ is a regular system of parameters in $\mathcal{O}_{Y,p}$. Because the residue field at $p$ is a finite separable extension of $k$, $dx_1,\dots,dx_n$ form a basis for the module $\Omega_{\mathcal{O}_{Y,p}/k}$ of K\"{a}hler differentials. By \cite[Corollary 6.2.11]{Liu06}, there exists $g \in A$ and an \'{e}tale $k$-algebra map $k[x_1,\dots,x_n] \to A[g^{-1}]$ such that the composition $k[x_1,\dots,x_n] \to A[g^{-1}] \into \mathcal{O}_{Y,p}$ maps the variable $x_i$ to the regular parameter $x_i$. By Proposition \ref{Rees algebra generators} and Lemma \ref{charts are affine}, the \'{e}tale map $k[x_1,\dots,x_n] \into A[g^{-1}]$ induces the following isomorphisms:
\begin{align*}
A_{\gamma}[g^{-1}]/(y_1-1) &= A[g^{-1}] [\{x_i,x_iT,\dots,x_iT^{w_i}\}_{i=1}^r] /(x_1T^{w_1}-1) \\
&= A[g^{-1}] \otimes_{k[x_1,\dots,x_n]}  k[x_1,\dots,x_n][\{x_i,x_iT,\dots,x_iT^{w_i}\}_{i=1}^r] /(x_1T^{w_1}-1)  \\
&=A[g^{-1}] \otimes_{k[x_1,\dots,x_n]} k[u,y_2,\dots,y_r,x_{r+1},\dots,x_n]
\end{align*}
Since \'{e}tale maps are stable under base change, we have the following \'{e}tale morphism over $k$:
$$
\Spec A_{\gamma}[g^{-1}]/(y_1-1) \to \Spec k[u,y_2,\dots,y_r,x_{r+1},\dots,x_n] =\mathbb{A}_k^n
$$
hence the distinguished open $\Spec A_{\gamma}[g^{-1}]/(y_1-1)$ of $W_1$ is regular and of dimension $n$. This implies that the pullback of $V(x_1,\dots,x_r)$ under the map $W_1 \to Y$ is covered by open affines that are regular.

The preimage of $\Spec A[x_1^{-1}]$ under the map $W_1 \to Y$ is $\Spec A_{\gamma}[x_1^{-1}]/(y_1-1)$. Note that $A_{\gamma}[x_1^{-1}] = A[x_1^{-1}][T]$, and that the map $A[x_1^{-1}] \to A[x_1^{-1}][T]/(x_1T^{w_1}-1)$ is standard {\'e}tale because $\frac{\partial}{\partial T} (x_1T^{w_1}-1) = w_1x_1T^{w_1-1}$ is invertible in $A[x_1^{-1}][T]/(x_1T^{w_1}-1)$ (note that $w_1$ is invertible because of the characteristic of $k$). Since $k \to A[x_1^{-1}]$ is smooth of relative dimension $n$, $\Spec A_{\gamma}[x_1^{-1}]/(y_1-1)$ is smooth over $k$ with dimension $n$.

Becuase $\Spec A[x_1^{-1}], \dots, \Spec A[x_r^{-1}]$, and $V(x_1,\dots,x_r)$ cover $Y$ and their pullbacks under $W_1\to Y$ are regular and of dimension $n$, it follows that $W_1$ is also regular and of dimension $n$. \qed

\subsection{The exceptional divisor} 
Note that $T$ is invertible in $A[T]/(x_1T^{w_1}-1)$, with inverse $T^{-1}=x_1T^{w_1-1}$, which is contained in $\mathcal{O}(W_1)=A^{pre}_{\gamma}[T^{-1}]/(x_1T^{w_1}-1)$. Let $u:=T^{-1} \in \mathcal{O}(W_1)$. By \cite[Lemma 4.4.3]{Que20}, $u$ is the equation of the exceptional divisor of the weighted blowings up.

\subsection{Weak and proper transforms} 
Let $I \subset A$ be an ideal, and let $\mathcal{I} \subset \mathcal{O}_Y$ be the associated coherent ideal. Without loss of generality, we compute the various transforms of the ideal $\mathcal{I}$ by the weighted blowings up restricted to the $x_1$-chart $W_1$. The total transform of $\mathcal{I}$ restricted to $W_1$ is given by the extension of the ideal $I$ under the map $A \to \mathcal{O}(W_1)$. The exceptional divisor is $u \in \mathcal{O}(W_1)$. Thus the proper transform of $\mathcal{I}$ restricted to $W_1$ is cut out by the ideal $(I \mathcal{O}(W_1) : u^{\infty})$. Note that the proper transform $V( I \mathcal{O}(W_1) : u^{\infty}) \subset W_1$ is indeed $\bmu_{w_1}$-invariant, because $I \mathcal{O}(W_1) \subset \mathcal{O}(W_1)$ lies in the zero graded piece and $u$ has degree one with respect to the $\mathbb{Z}/w_1$-grading of $\mathcal{O}(W_1)$. Thus the proper transform of $V(\mathcal{I}) \subset Y$ restricted to the open substack $[W_1 / \bmu_{w_1}]$ is
$$
[ V( I \mathcal{O}(W_1) : u^{\infty}) / \bmu_{w_1}] \subset [W_1 / \bmu_{w_1}]
$$
Now assume that $\gamma=(x_1^{1/w_1},\dots,x_r^{1/w_r})$ is the reduced center associated to $\mathcal{I}$, and that $\maxinv \mathcal{I} = (a_1,\dots,a_r)$. Then the ideal of the weak transform of $\mathcal{I}$ restricted to $W_1$ is $(I \mathcal{O}(W_1) : u^{a_1w_1})$. See \cite[5.2.1]{ATW20}.

\subsection{Algorithm for weighted blowing ups} \
\Suppressnumber
\begin{lstlisting}[caption={weightedBlowUp},label={weightedBlowUp}]
Let $k$ be a field of characteristic zero.
$\text{\bf{Input:}}$ Ideals $I,I_Y\subset k[x_1,\dots,x_N]$, polynomials $f_1,\dots,f_r \in k[x_1,\dots,x_N]$, and positive integers $w_1,\dots,w_r \in \mathbb{N}$, such that $Y= \Spec A$ is a smooth irreducible variety over $k$, where $A=k[x_1,\dots,x_N]/I_Y$, and $f_1,\dots,f_r$ are regular parameters on $Y$. Let $\mathcal{I} \subset \OY$ be the coherent ideal on $Y$ associated to $I$. 
$\text{\bf{Compute:}}$ The charts $W_1,\dots,W_r$ of the weighted blowings up of $Y$ along the center $(f_1^{1/w_1},\dots,f_r^{1/w_r})$, where each $W_i$ is the @$f_i$-chart@, as well as the proper transform of $\mathcal{I}$. 
$\text{\bf{Output:}}$ A list $\text{CHARTS}$ of size $r$ such that for each $i$, $\text{CHARTS}[i] = (R_i, I_{W_i}, I_i, \phi_i)$, where $R_i$ is a polynomial ring over $k$, $I_{W_i}, I_i \subset R_i$ are ideals, and $\phi_i: k[x_1,\dots,x_N] \to R_i$ is a @$k$-algebra@ map, such that the coordinate ring of the @$f_i$-chart@ $W_i$ is $R_i/I_{W_i}$, the ideal of the proper transform of $\mathcal{I}$ on the @$f_i$-chart@ is the image of $I_i$ in $R_i/I_{W_i}$, and $\phi_i$ induces a map $A \to R_i/I_{W_i}$, which induces the weighted blowings up morphism on the @$f_i$-chart@ $W_i \to Y$.@\Reactivatenumber@ 
Introducing new variables $y_1,\dots,y_r,T,u$, define the following @$k[x_1,\dots,x_N]$-algebra@ map @$$\phi: k[x_1,\dots,x_N,y_1,\dots,y_r, u] \to k[x_1,\dots,x_N,T,u]$$@given by $y_i \mapsto f_iT^{w_i}$ and $u \mapsto u$.@\footnote{See {\sc Singular}'s procedure \texttt{map} \cite[4.10]{DGPS19}}@
Denote the image of $I_Y$ in $k[x_1,\dots,x_N,T,u]$ also by $I_Y$. Let@\footnote{See {\sc Singular}'s procedure \texttt{preimage} \cite[5.1.111]{DGPS19}}@  @$$I_{\phi} = \phi^{-1}\big((I_Y,Tu-1)\big) \subset k[x_1,\dots,x_N,y_1,\dots,y_r, u]$$@
Observe that @$$\frac{k[x_1,\dots,x_N,T,u]}{(I_Y,Tu-1)}=\bigg(\frac{k[x_1,\dots,x_N]}{I_Y}\bigg)[T,u]/(Tu-1) =A[T,T^{-1}]$$@where $u=T^{-1}$. Thus the inclusion @$$k[x_1,\dots,x_N,y_1,\dots,y_r, u]/I_{\phi} \into A[T,T^{-1}]$$@induced by $\phi$ identifies $k[x_1,\dots,x_N,y_1,\dots,y_r, u]/I_{\phi}$ with the @$A$-subalgebra@ of $A[T,T^{-1}]$ generated by $f_1T^{w_1},\dots, f_rT^{w_r},T^{-1}$. This means that @$$k[x_1,\dots,x_N,y_1,\dots,y_r, u]/I_{\phi}=A^{pre}_{\gamma}[T^{-1}]\subset A[T,T^{-1}]$$@
$\text{\bf{for}}$ each $i=1,\dots,r$ $\text{\bf{do}}$
  We compute the @$x_i$-chart@. By Corollary @\ref{compute coordinate ring}@, the $f_i$-chart $W_i$ is given by @$$\mathcal{O}(W_i) = A^{pre}[T^{-1}]/(x_iT^{W_i}-1)=\frac{k[x_1,\dots,x_N,y_1,\dots,y_r, u]}{(I_{\phi},f_iT^{w_i}-1)}$$@
  We can eliminate the variable $y_i$, which is equal to one in $\mathcal{O}(W_i)$, as well as maybe some of the variables among $x_1,\dots,x_N$.@\footnote{See {\sc Singular}'s procedure \texttt{elim} from the library \texttt{elim.lib} \cite[D.4.5.3]{DGPS19}.}@ Let $E \subset \{1,\dots,N\}$ be the indices such that the variables $\{ x_j \}_{j \in E}$ can be eliminated. This means there are polynomials $\{g_j\}_{j \in E} \subset k[x_1,\dots,x_N,y_1,\dots,y_r, u]$ such that for each $j \in E$, $z_j - g_j \in (I_{\phi},f_iT^{w_i}-1)$ and $g_j$ does not involve the variables $\{ x_j \}_{j \in E}\cup \{y_i\}$. Let @$$\text{elim}_i:  k[x_1,\dots,x_N,y_1,\dots,y_r, u] \to  k[\{x_{\ell}\}_{\ell \not\in E},\{ y_j\}_{j \neq i}, u]$$@be the @$k[\{x_{\ell}\}_{\ell \not\in E},\{y_j\}_{j\neq i}, u]$-algebra@ map sending $z_j \mapsto g_j$ for each $j \in E$ and $y_i \mapsto 1$, and let the ideal $I_{W_i}$ be the extension of $(I_{\phi},f_iT^{w_i}-1)$ under $\text{elim}_i$. Then $\text{elim}_i$ induces the @$k$-algebra@ isomorphism @$$\frac{k[x_1,\dots,x_N,y_1,\dots,y_r, u]}{(I_{\phi},f_iT^{w_i}-1)} \xrightarrow{\cong} \frac{k[\{z_{\ell}\}_{\ell \not\in E},\{y_j\}_{j \neq i}, u]}{I_{W_i}}$$@Thus $k[\{z_{\ell}\}_{\ell \not\in E},\{y_j\}_{j \neq i} u]/I_{W_i}$ is the coordinate ring of $W_i$.
  Let $\phi_i$ be the composition @$$k[x_1,\dots,x_N] \into k[x_1,\dots,x_N,y_1,\dots,y_r, u] \xrightarrow{\text{elim}_i}k[\{x_{\ell}\}_{\ell \not\in E},\{y_j\}_{j \neq i}, u]$$@and let $\text{Tot}_i \subset k[\{x_{\ell}\}_{\ell \not\in E},\{y_j\}_{j \neq i}, u]$ be the image of $I$ under $\phi_i$. $\phi_i$ induces the map @$$k[x_1,\dots,x_N]/I_Y\to k[\{x_{\ell}\}_{\ell \not\in E},\{y_j\}_{j \neq i}, u]/I_{W_i}$$@which is the map of coordinate rings induced by $W_i \to Y$.
  The image of $\text{Tot}_i$ in $k[\{x_{\ell}\}_{\ell \not\in E},\{y_j\}_{j \neq i}, u]/I_{W_i}$ is the ideal $I\mathcal{O}(W_i)$ of the total transform of $I$ on the @$f_i$-chart@, while $u$ is the equation of the exceptional divisor. Let @$$I_i = (\text{Tot}_i + I_{W_i}: u^{\infty}) \subset k[\{z_{\ell}\}_{\ell \not\in E},\{y_j\}_{j \neq i}, u]$$@The image of $I_i$ in $k[\{z_{\ell}\}_{\ell \not\in E},\{y_j\}_{j \neq i}, u]/I_{W_i}$ is the ideal $(I\mathcal{O}(W_i):u^{\infty})$ of the proper transform of $\mathcal{I}$ on the @$f_i$-chart@. 
  $\text{CHARTS}[i]=(k[\{x_{\ell}\}_{\ell \not\in E},\{y_j\}_{j \neq i},u], I_{W_i}, I_i, \phi_i)$
  $\text{\bf{Remark:}}$ The @$\bmu_{w_i}$-action@ on $W_i$ is recovered from the presentation @$$k[\{x_{\ell}\}_{\ell \not\in E},\{y_j\}_{j \neq i},u]/I_{W_i}$$@ of the coordinate ring of $W_i$ by placing $\{x_{\ell}\}_{\ell \not\in E}$ in degree zero, each $y_j$ in degree $-w_j$, and $u$ in degree one. Hence the stacky structure of the charts is inherent in the presentation of their coordinate rings.
$\text{\bf{return}}$ $\text{CHARTS}$
\end{lstlisting}

\subsection{Example for weightedBlowUp}
We illustrate Algorithm \nameref{weightedBlowUp} with the weighted blowings up of $Y = \Spec k[x,y]$ along the center $(x^{1/3},y^{1/2})$, as well as the proper transform of the ideal $(x^5+x^3y^3+y^8)$. 

First, define the $k[x,y]$-algebra map 
\begin{align*}
\phi: k[x,y,y_1,y_2,u] &\to k[x,y,T,u] \\
y_1 &\mapsto xT^3 \\
y_2 &\mapsto yT^2 \\
u &\mapsto u 
\end{align*}
Let $I_{\phi} = \phi^{-1}\big( (Tu-1) \big)$ be the preimage of the ideal $(Tu-1)$. By {\sc Singular}'s procedure \texttt{preimage} \cite[5.1.111]{DGPS19}, we obtain
$$
I_{\phi} =(x-y_1u^3, \ y-y_2u^2, \ y^3y^2_1-x^2y_2^3, \ yy_1u-xy_2, \ xy_2^2u-y^2y_1)
$$
Note that $\phi$ identifies $k[x,y,y_1,y_2,u]/I_{\phi}$ with the $k[x,y]$-subalgebra of $k[x,y,T,u]/(Tu-1)=k[x,y][T,T^{-1}]$ generated by $xT^3,yT^2,T^{-1}$. 

The coordinate ring of the $x$-chart $W_1$ is given by 
$$
\frac{k[x,y,y_1,y_2,u]}{(I_{\phi}, y_1-1)}
$$
By {\sc Singular}'s procedure \texttt{elim} from the library \texttt{elim.lib} \cite[D.4.5.3]{DGPS19}, we can eliminate the variables $x,y,y_1$ through the equations $x-u^3, y-y_2u^2, y_1-1 \in (I_{\phi}, y_1-1)$. Consider the substitution $k[y_2,u]$-algebra map
\begin{align*}
\text{elim}_1: k[x,y,y_1,y_2,u] &\to k[y_2,u] \\
x &\mapsto u^3 \\
y &\mapsto y_2u^2 \\
y_1 &\mapsto 1
\end{align*}
Let the ideal $I_{W_1} \subset k[y_2,u]$ be the extension of $(I_{\phi}, y_1-1)$ under the map $\text{elim}_1$. Note that $I_{W_1}=0$. Thus the $x$-chart is $W_1=\Spec k[y_2,u]$, and the morphism $W_1 \to Y$ is given by the following $k$-algebra composition
$$
k[x,y] \into k[x,y,y_1,y_2,u] \xrightarrow{\text{elim}_1} k[y_2,u]
$$
which sends $x\mapsto u^3$ and $y\mapsto y_2u^2$. The total transform of $x^5+x^3y^3+y^8$ on the $x$-chart is:
$$
(u^3)^5+(u^3)^3(y_2u^2)^3+(y_2u^2)^8=u^{15}(1+y_2^3+uy_2^8)
$$
where $u$ is the equation of the exceptional divisor. Thus $1+y_2^3+uy_2^8$ is the proper transform of $x^5+x^3y^3+y^8$ on the $x$-chart. The $\mu_3$-action on the $x$-chart is given by placing $y_2$ in degree $-2\equiv 1$ and $u$ in degree $1$. 

Similarly, we obtain that the $y$-chart is $W_2=\Spec k[y_1,u]$, where the morphism $W_2 \to Y$ is given by the $k$-algebra map $k[x,y] \to k[y_1,u]$ sending $x\mapsto y_1u^3$ and $y\mapsto u^2$, and that the proper transform of $x^5+x^3y^3+y^8$ on the $y$-chart is $y_1^5+y_1^3+u$. Also, the $\mu_2$-action on the $y$-chart is given by placing $y_1$ in degree $-3\equiv 1$ and $u$ in degree $1$.


\section{Weighted resolution algorithm}
\label{sec: resweighted}

\subsection{Representing the input}
Before we explain the algorithm for embedded resolution, we need to discuss how we want to represent varieties, or more generally Deligne-Mumford stacks, as computer code. The nature of the embedded resolution functor $F_{er}$ of \cite[Theorem 1.1.1]{ATW20} will guide us how to do so. 

Let $Y$ be a smooth Deligne-Mumford stack over a field $k$ of characteristic zero, and $X \subset Y$ a closed substack. How should we tell the computer what the pair $(X \subset Y)$ is? Well, suppose $(\tilde{X} \subset \tilde{Y})$ is an {\'e}tale cover of $(X \subset Y)$ by schemes, that is, $\tilde{X}$ and $\tilde{Y}$ are varieties with an {\'e}tale surjection $\tilde{Y} \to Y$ such that $\tilde{X}=X \times_Y \tilde{Y}$. By functoriality with respect to smooth surjections, $F_{er}(\tilde{X} \subset \tilde{Y})$ is an {\'e}tale cover of $F_{er}(X \subset Y)$, and by \cite[Remark 2.1.1]{ATW20}, $F_{er}(X \subset Y)$ is obtained by {\'e}tale descent from $F_{er}(\tilde{X} \subset \tilde{Y})$. Hence functoriality with respect to smooth surjections implies that the embedded resolution of $(X \subset Y)$ is fully determined by its {\'e}tale cover $(\tilde{X} \subset \tilde{Y})$. Thus we will represent the data of the pair $(X \subset Y)$ by its {\'e}tale cover $(\tilde{X} \subset \tilde{Y})$ by schemes. 

Now the question shifts to how we will represent the pair $(X \subset Y)$ if $X$ and $Y$ are varieties. Suppose $Y= \cup_i Y_i$ is a cover of $Y$ by finitely many open affines $Y_i$ and ideals $I_i \subset \mathcal{O}(Y_i)$ such that $I_i$ is the defining ideal of the closed subvariety $X|_{Y_i}$ in $Y_i$. We can input the affine varieties $Y_i$ and the ideals $I_i$ into a computer as long as we know some presentation of each of the coordinate rings $\mathcal{O}(Y_i)$. Let $\tilde{Y} = \sqcup_i Y_i$, and let $\tilde{X} \subset \tilde{Y}$ be the closed subvariety defined by the ideals $I_i \subset \mathcal{O}(Y_i)$. $(\tilde{X} \subset \tilde{Y})$ is an {\'e}tale cover of $(X \subset Y)$ by schemes, and hence we represent the data of $(X \subset Y)$ by its {\'e}tale cover $(\tilde{X} \subset \tilde{Y})$. Since the data of $(\tilde{X} \subset \tilde{Y})$ is equivalent to the data of the coordinate rings $\mathcal{O}(Y_i)$ and the ideals $I_i \subset \mathcal{O}(Y_i)$, there is no problem inputting $(\tilde{X} \subset \tilde{Y})$ into a computer.

We do note that in passing to {\'e}tale covers we lose information about gluing, as can be seen when passing from $(X \subset Y)$ to $(\tilde{X} \subset \tilde{Y})$. However, we will not concern ourselves with this matter, since the embedded resolution functor $F_{er}$ operates {\'e}tale-locally. 

\subsection{Examples of representing inputs}

For example, the pair $(V(zx^2-y^3) \subset \mathbb{P}^2_{k})$ can be inputted into the weighted resolution algorithm by the data of ideals 
\begin{align*}
(z-y^3) &\subset k[z,y] \\
(zx^2-1) &\subset k[x,z] \\
(x^2-y^3) &\subset k[x,y]
\end{align*}
which are the ideals defining the projective subvariety $V(zx^2-y^3)$ on the standard $x$, $y$, and $z$-charts of $\mathbb{P}^2_{k}=\Proj k[x,y,z]$, respectively. 

For a stacky example, consider the pair $( [V(x^3-y^2)/\bmu_{5}]  \subset [\mathbb{A}^2_k/\bmu_{5}] )$, where the $\bmu_{5}$-action on $\mathbb{A}^2_k=\Spec k[x,y]$ is given by placing $x$ in degree $2$ and $y$ in degree $3$. Instead of inputing the pair $( [V(x^3-y^2)/\bmu_{5}]  \subset [\mathbb{A}^2_k/\bmu_{5}] )$ into weighted resolution, we input its {\'e}tale cover $(V(x^3-y^2) \subset \mathbb{A}^2_k)$.

\subsection{Algorithm for weighted resolution of singularities} \
\Suppressnumber
\begin{lstlisting}[caption={weightedResolution},label={weightedResolution}]
Let $k$ be a field of characteristic zero.
$\text{\bf{Input:}}$ $Y$ a smooth irreducible variety over $k$ and $X \subset Y$ a pure codimension reduced closed subvariety.
$\text{\bf{Compute:}}$ The embedded resolution of singularities of the pair $(X \subset Y)$. That is, we compute the stabilized functor $F^{\circ \infty}_{er}( X \subset Y)$.@\footnote{See \cite[Theorem 1.1.1]{ATW20}}@ 
$\text{\bf{Output:}}$ A list called $\Res$ such that
@{\bu}@$n=\text{size}(\Res)$ is the smallest integer such that $(X_n \subset Y_n):=F^{\circ n}_{er}( X \subset Y)$ has $X_n$ smooth. For example, if $X$ is already smooth, then $\Res$ is the empty list.
@{\bu}@For $i=1,\dots,n$, $\Res[i]=(X_i \subset Y_i)$, where $(X_i \subset Y_i):=F^{\circ i}_{er}( X \subset Y)$. 
$\text{\bf{Remark:}}$ We will not distinguish between a pair of Deligne-Mumford stacks and its @{\'e}tale@ cover by schemes. But keep in mind, $\Res[i]$ is really the data of an @{\'e}tale@ cover of the pair $(X_i \subset Y_i)$ of Deligne-Mumford stacks by schemes. @\Reactivatenumber@ 
Let $c=\dim Y-\dim X$, which is the codimension of $X$ in $Y$.
Let $\mathcal{I} \subset \mathcal{O}_{Y}$ be the coherent ideal on $Y$ cutting out $X$.
Initialization
  @{\bu}@$i=0$
  @{\bu}@$\Res=\emptyset$
  @{\bu}@$\text{SingularitiesNotResolved}=\T$
Let $(X_{curr}\subset Y_{curr}) = (X \subset Y)$, and let $\mathcal{I}_{curr}=\mathcal{I}$. 
$\text{\bf{while}}$ $\text{SingularitiesNotResolved}$ $\text{\bf{do}}$
  $i=i+1$
  Compute $\maxinv \mathcal{I}_{curr}$ and the reduced center $\bar{J}$ on $Y_{curr}$ associated to $\mathcal{I}_{curr}$ using Algorithm @\nameref{prepareCenter}@.
  $\text{\bf{if}}$ $\maxinv \mathcal{I}_{curr} = (1,\dots,1)$, the constant sequence of length $c$ $\text{\bf{do}}$
    $X_{curr}=V(\mathcal{I}_{curr})$ is smooth!
    $\text{SingularitiesNotResolved}=\F$
  $\text{\bf{else}}$
    $X_{curr}=V(\mathcal{I}_{curr})$ is not smooth.
    Let $Y_{new}$ be the weighted blowings up of the reduced center $\bar{J}$ on $Y_{curr}$, and let $X_{new}$ be the proper transform of $X_{curr}$.
    $\Res[i]=(X_{new}\subset Y_{new})$
    Set $Y_{curr}=Y_{new}$, $X_{curr}=X_{new}$.
$\text{\bf{return}}$ $\Res$
\end{lstlisting}

\subsection{Example for weightedResolution}
We illustrate Algorithm \nameref{weightedResolution} with the weighted resolution of the Whitney umbrella $X=V(x^2-zy^2)$ in $Y=\Spec k[x,y,z]$. 

By Algorithm \nameref{prepareCenter}, we have that $\maxinv X =(2,3,3)$ and that the reduced center associated to $(x^2-zy^2)$ is $(x^{1/3},y^{1/2},z^{1/2})$. We need to take the weighted blowings up of $Y$ along $(x^{1/3},y^{1/2},z^{1/2})$. 
\begin{itemize}
\item In the $x$-chart $W_1=\Spec k[y_2,y_3,u]$, the equation of $X$ becomes 
$$
u^6(1-y_2y_3^2)
$$ 
with proper transform $X'_x=V(1-y_2y_3^2)$ smooth.
\item In the $y$-chart $W_2 =\Spec k[y_1,y_3,u]$, the equation of $X$ becomes 
$$
u^6(y_1^2-y_3)
$$ 
with proper transform $X'_y=V(y_1^2-y_3)$ smooth.
\item In the $z$-chart $W_3=\Spec k[y_1,y_2,u]$, the equation of $X$ becomes
$$
u^6(y_1^2-y_2^2)
$$
with proper transform $X'_z=V(y_1^2-y_2^2)$. We have that $\maxinv X'_z = (2,2)$, and that the reduced center associated to $(y_1^2-y_2^2)$ is $(y_1,y_2)$. The weighted blowings up of $W_3$ along the center $(y_1,y_2)$ resolves the singularities.
\end{itemize}


\section{Comparisons with Villamayor's resolution algorithm}
\label{sec: compare}

In this section, we will compare the weighted resolution algorithm with Villamayor's resolution algorithm, in particular with respect to efficiency. Both of these algorithms are implemented in {\sc Singular}: the library \texttt{resweighted.lib} implements weighted resolution \cite{AFL20}, and the library \texttt{resolve.lib} implements Villamayor's resolution algorithm \cite{FP19}. Through these computer implementations, we will compare the efficiencies of the two resolution algorithms.

\subsection{Two criterions of efficiency}
Both algorithms represent varieties by their affine chart coverings and compute (weighted or non-weighted) blowing ups locally on these affine charts, leading us towards the following two criterions of efficiency: the number of times a (weighted or non-weighted) blowings up is performed and the number of charts covering the final resolution. For Villamayor's algorithm, we obtain the number of times a blowings up is performed by counting the number of times the procedure \texttt{blowUpBo} is called within the procedure \texttt{resolve}, both of which are procedures in \texttt{resolve.lib}. Similarly, for weighted resolution, we obtain the number of times a weighted blowings up is performed by counting the number of times the procedure \nameref{weightedBlowUp} is called within the procedure \nameref{weightedResolution}, both of which are procedures in \texttt{resweighted.lib}. Note that the number of blowing ups and the number of final charts are related, since more blowing ups contributes to increasing the number of charts. 

\subsection{Experimental results}
We tabulate the results of our computer trials in the table below. The column under ``Input pair" represents the input data: for each polynomial $f$ in this column, we input into both resolution algorithms the pair $(V(f) \subset \Spec k[x_1,\dots,x_n])$, where $x_1,\dots,x_n$ are the variables involved in $f$. For example, $x^2-y^2z$ represents the pair $(V(x^2-y^2z) \subset \Spec k[x,y,z])$. The columns under ``Weighted resolution" and ``Villamayor's algorithm" contain the measurements of the two criterions of efficiency for each corresponding input pair, where the ordered pairs in these columns follow the rule (\# of blowing ups, \# of final charts). 

\setlength\extrarowheight{5pt}
\begin{longtable}{|l|l|l|}
\hline
\textbf{Input pair}                       & \textbf{Weighted resolution}              & \textbf{Villamayor's algorithm} \\ \hline
$x^2-y_1y_2y_3$ & (5,10)  & (166,198)          \\ \hline
$x^2+y^2+z^3t^3$ &  (5,8)  &    (27,35)           \\ \hline
$x^5+x^3y^3+y^7$ & (1,2) & (5,2) \\ \hline
$x^5+x^3y^3+y^9$ & (3,3) & (3,2) \\ \hline
$x^2+y^2z^3-z^4$ & (4,5) & (9,10) \\ \hline
$x^2+y^2z-z^2$ & (1,3) & (2,3) \\ \hline
$x^3y+xz^3+y^3z+z^3+7z^2+5z$ & (1,3)& (4,5)\\ \hline
$x^2+y^3+z^5$ & (1,3) & (61,65) \\ \hline
$x^2-x^3+y^2+y^4+z^3-z^4$ & (1,3) & (4,5) \\ \hline
$x^3-y(1-z^2)^2$ & (4,4) & (50,49) \\ \hline
$x^4+z^3-yz^2$ & (4,4) & (50,49) \\ \hline
$x^2+z^2+y^3(y-1)^3$ & (1,3) & (8,10) \\ \hline
$xyz+yz+2z^5$ & (9,7) & (6,10) \\ \hline
$x^2+y^4+y^3z^2$ & (4,5) & (9,10) \\ \hline
$x^2+y^2z^3$ & (8,5) & (24,23) \\ \hline
$xyz$ & (4,6) & (7,12) \\ \hline
$x^2+y^2z+z^3$ & (1,3) & (22,24) \\ \hline
$z^50-xy$ & (1,3) & (25,50) \\ \hline
\end{longtable}

\subsection{Conclusion}
For each input pair, we see that weighted resolution is either on par or vastly superior than Villamayor's algorithm in both criterions of efficiency. For example, the weighted resolution of the pair $(V(x^2+y^3+z^5) \subset \Spec k[x,y,z])$ requires just one weighted blowings up, with three final charts, while the resolution of the pair given by Villamayor's algorithm requires $61$ blowing ups, with $65$ final charts. From these computer trials, we therefore conclude that weighted resolution is more efficient than Villamayor's algorithm. Indeed, weighted resolution is remarkably efficient!

\appendix

\section{Zariski closure and Saturation}
\label{sec: saturation}

In this appendix, we present a scheme-theoretic treatment of Zariski closure and saturation. See \cite[\S 4.4]{CLO15} for a treatment of Zariski closure and saturation in the context of classical varieties.

\subsection{Zariski closure}

\begin{definition}[Scheme-Theoretic Image and Zariski closure, see {\cite[\S 8.3]{Vak17}}] Let $f: X \to Y$ be a morphism of schemes. The scheme-theoretic image of $f$ is the smallest closed subscheme $Z \subset Y$ through which $f$ factors through. The scheme-theoretic image $Z$ is cut out by the sum of all quasicoherent ideals contained in $\ker(\mathcal{O}_Y \to f_*\mathcal{O}_X)$. If $f: X \to Y$ is a locally closed embedding, then we call the scheme-theoretic image of $f$ the scheme-theoretic closure or the Zariski closure, and we denote the Zariski closure of $X \into Y$ by $\overline{X}$. 
\end{definition}

If $B \to A$ is a ring map, then the kernel of $B \to A$ cuts out the scheme-theoretic image of the morphism $\Spec A \to \Spec B$. 

The following lemma and corollary are general results that will be needed for Section \ref{sec: sat}.

\begin{lemma} 
\label{kernel intersection}
Let $f: X \to Y$ be a morphism of schemes and let $X = \bigcup_{i\in I} X_i$ be a covering of $X$ by open subschemes $X_i$. Then
$$
\bigcap_{i \in I} \ker\big( \mathcal{O}_Y \to (f|_{X_i})_* \mathcal{O}_{X_i} \big) = \ker ( \mathcal{O}_Y \to f_* \mathcal{O}_X)
$$
\end{lemma}
\proof This follows directly from the definition of sheaves. \qed

\begin{corollary} 
\label{glue Zariski closures}
Let $f: X \to Y$ be a morphism of schemes and let $X = \bigcup_{i=1}^r X_i$ be an open covering of $X$ by finitely many open subschemes $X_i$. If $\mathcal{I} \subset \OY$ cuts out the scheme-theoretic image of $f: X \to Y$ and $\mathcal{I}_i \subset \OY$ cuts out the scheme-theoretic image of $f|_{X_i} : X_i \into X \to Y$, then $\mathcal{I} = \bigcap_{i=1}^r \mathcal{I}_i$. 
\end{corollary}
\proof 
By Lemma \ref{kernel intersection}, 
$$
\ker ( \mathcal{O}_Y \to f_* \mathcal{O}_X) = \bigcap_{i=1}^r \ker\big( \mathcal{O}_Y \to (f|_{X_i})_* \mathcal{O}_{X_i} \big) = \bigcap_{i=1}^r \mathcal{I}_i
$$
Thus $\ker ( \mathcal{O}_Y \to f_* \mathcal{O}_X)$ is quasicoherent. Hence we get $\mathcal{I} = \ker ( \mathcal{O}_Y \to f_* \mathcal{O}_X)=\bigcap_{i=1}^r \mathcal{I}_i$. \qed

\subsection{Saturation}
\label{sec: sat}

\begin{definition}[Saturation] Let $A$ be a ring and $I,J \subset A$ ideals. Define the saturation $(I: J^{\infty})$ of $I$ with respect to $J$ to be the following ideal
$$
(I: J^{\infty}) := \bigcup_{n=1}^{\infty} (I:J^n)
$$
If $h\in A$, then define $(I: h^{\infty}) := (I: (h)^{\infty})$.\footnote{Use {\sc Singular}'s procedure \texttt{sat} from the \texttt{elim.lib} library to compute saturation \cite[D.4.5.7]{DGPS19}}
\end{definition}

Note that $(I: h^{\infty}) = \ker (A \to (A/I)[h^{-1}])$. For the following proposition, let $o.e.$ and $c.e.$ stand for open and closed embedding, respectively. Recall that morphisms that factor as an open embedding followed by a closed embedding are locally closed embeddings \cite[Exercise 8.1.M]{Vak17}.

\begin{proposition}[Geometric Interpretation of Saturation] Let $A$ be a ring and $I,J \subset A$ ideals, where $J$ is a finitely generated ideal. Then the Zariski closure of $V(I) \setminus V(J) \xhookrightarrow{o.e.} V(I) \xhookrightarrow{c.e.} \Spec A$ is cut out by $(I:J^{\infty}) \subset A$. 
\end{proposition}
\proof First consider the case that $J$ is a principal ideal generated by $h \in A$. The open embedding $V(I) \setminus V(h) \into V(I)$ is precisely the inclusion of the distinguished affine open $\Spec (A/I)[h^{-1}] \into \Spec A/I$. Thus the Zariski closure of $V(I) \setminus V(h)$ in $\Spec A$ is cut out by $\ker(A \to (A/I)[h^{-1}]) = (I: h^{\infty})$. 

Now suppose that $J$ is not necessarily principal and generated by $h_1,\dots,h_r \in A$. Then 
$$
V(I) \setminus V(J) = V(I) \setminus \bigcap_{i=1}^r V(h_r) = \bigcup_{i=1}^r V(I) \setminus V(h_r)
$$
so that $V(I) \setminus V(J)$ is covered by finitely many open affines $V(I) \setminus V(h_r)$. By Corollary \ref{glue Zariski closures}, the Zariski closure of $V(I) \setminus V(J)$ in $\Spec A$ is cut out by $\bigcap_{i=1}^r (I: h_i^{\infty})$. We claim that
$$
\bigcap_{i=1}^r (I: h_i^{\infty}) = (I: (h_1, \dots, h_r)^{\infty})
$$
To see this, first let $f \in (I: (h_1, \dots, h_r)^{\infty})$. Then $f (h_1, \dots, h_r)^n \subset I$ for some $n$, so that $f \in (I: h_i^n)$ for every $i$. Conversely, suppose $f \in \bigcap_{i=1}^r (I: h_i^{\infty})$. Then for sufficiently large $N$, we have $fh_i^N \in I$ for each $i$, so that $f \in (I: (h_1, \dots, h_r)^{rN})$. Thus we are done. \qed

\begin{proposition}[Gluing Ideals] 
\label{gluing ideals}
Let $A$ be a ring and $h_1, \dots, h_r \in A$ be elements generating the unit ideal, i.e. $(h_1, \dots, h_r)=A$. Let $I \subset A$ and $I_i \subset A$ for $i =1, \dots,r$ be ideals such that $I$ and $I_i$ have the same extension in $A[h_i^{-1}]$ for each $i$. Then
$$
I = \bigcap_{i=1}^r ( I_i : h_i^{\infty})
$$
\end{proposition}
\proof First we give a geometric picture. Since $V(I)$ and $V(I_i)$ agree on $D(h_i)$, we have $V(I) \subset \cup V(I_i)$. There may be components of $V(I_i)$ that lie entirely in $V(h_i)$ that do not appear in $V(I)$. So we remove these extraneous components by taking the Zariski closure of $V(I_i) \setminus V(h_i)$. Thus $V(I) = \bigcup_{i=1}^r \overline{V(I_i) \setminus V(h_i)}$. Scheme-theoretically, this translates to precisely the algebraic statement. 

Now the algebraic proof. If $f \in I$, then because $I$ and $I_i$ have the same extension in $A[h_i^{-1}]$, there is $N_i$ such that $f h_i^{N_i} \in I_i$, hence $f \in (I_i : h_i^{\infty})$.

Conversely, assume that $f \in \bigcap_{i=1}^r ( I_i : h_i^{\infty})$. Then for some sufficiently high $N$, we have $fh_i^N \in I_i$ for each $i$. Since $I$ and $I_i$ have the same extension in $A[h_i^{-1}]$, by taking sufficiently larger $N$, we have $fh_i^N \in I$ for each $i$. Because $(h_1, \dots, h_r)=A$, there exists $a_i \in A$ such that $\sum_{i=1}^r a_i h_i^N = 1$ (``partition of unity" trick). Thus
$$
f = \sum_{i=1}^r a_i fh_i^N \in I
$$
\qed

\subsection{Computational remark}
\label{subsec: compute remark}
The procedure \texttt{sat} in the {\sc Singular} library \texttt{elim.lib} computes the saturation $(I:J^{\infty})$ for any two ideals $I,J \subset k[x_1,\dots,x_N]$ in a polynomial ring \cite[D.4.5.7]{DGPS19}. We can use the procedure \texttt{sat} to compute saturation in quotients of polynomial rings as well, but we will first need the following general observation.

Let $R$ be a ring and $I,J,I_Y \subset R$ ideals such that $I_Y \subset I$. Let the superscript $e$ denote extensions of ideals along the quotient map $R\to R/I_Y$. Observe that $(I: J)^e = (I^e : J^e)$ as ideals of $R/I_Y$ (this observation uses the inclusion $I_Y \subset I$ in an essential way). Thus $(I: J^{\infty})^e = (I^e : (J^e)^{\infty})$. 

Now set $R=k[x_1,\dots,x_N]$. Every ideal of $R/I_Y$ is an extension of an ideal in $R$, hence by the above observation, we can use the procedure \texttt{sat} to compute saturation of ideals in the finitely generated $k$-algebra $R/I_Y$. 

We also make a remark on computing intersection of ideals in finitely generated $k$-algebras, because Algorithm \nameref{Diff} relies on Proposition \ref{gluing ideals}. Let $R$ be a ring and $I_Y, I_1,I_2 \subset R$ ideals such that $I_i \supset I_Y$ for each $i$. Observe that $(I_1 \cap I_2)^e = I_1^e \cap I_2^e$, where the superscript $e$ denote extensions of ideals along the quotient map $R\to R/I_Y$. Now set $R=k[x_1,\dots,x_N]$. Since every ideal of $R/I_Y$ is an extension of an ideal in $R$, we can use the {\sc Singular} procedure \texttt{intersect} \cite[5.1.60]{DGPS19}, which computes intersections of ideals in a polynomial ring, to compute intersection of ideals in the finitely generated $k$-algebra $R/I_Y$.

\section{Helper algorithms}
\label{sec: helper}

\subsection{Algorithm for determining whether an equation cuts out a smooth hypersurface}
\label{subsec: smooth hypersurface}
Let $Y$ be a variety. We say that a subvariety $X \subset Y$ is a hypersurface of $Y$ if $X$ is pure codimension one in $Y$, that is, all irreducible components of $X$ are codimension one in $Y$. The following algorithm determines whether or not an equation cuts a smooth hypersurface out of a smooth affine variety.
\Suppressnumber
\begin{lstlisting}[caption={isSmoothHypersurface},label={isSmoothHypersurface}]
$\text{\bf{Input:}}$ $f_1,\dots,f_r, h \in k[x_1,\dots,x_N]$ such that $Y=\Spec A$ is a smooth variety over $k$, where $A = k[x_1,\dots,x_N]/(f_1,\dots,f_r)$. Let $V(h) = \Spec A/h$ be the subvariety cut out by $h$.
$\text{\bf{Output:}}$ $\T$ if $V(h)$ is a smooth hypersurface, that is, a smooth subvariety of pure codimension one in $Y$, and $\F$ otherwise. @\Reactivatenumber@ 
By Krull's principal ideal theorem, $V(h)$ is pure codimension one in $Y$ iff $h$ does not vanish on any irreducible component of $Y$. 
Let $\frak{p}_i \subset k[x_1,\dots,x_N]$ be the associated primes of the ideal $(f_1,\dots,f_r) \subset k[x_1,\dots,x_N]$.@\footnote{See {\sc Singular}'s primary decomposition library \cite{PDSL19}}@
Because $\Spec A$ is a smooth variety, the associated primes $\{\frak{p}_i\}_i$ are the minimal primes of $A$ and $\{ \Spec A/\frak{p}_i \}_i$ are the irreducible components of $\Spec A$.
$\text{\bf{if}}$ there is an $i$ such that $h \in \frak{p}_i$ $\text{\bf{do}}$
  In this case, the subvariety $V(h)$ contains the component $\Spec A/\frak{p}_i$.
  $\text{\bf{return}}$ $\F$
Since we are here, $V(h)$ is indeed a hypersurface of $Y$. 
Now we verify whether $V(h)$ is smooth or not.
Let $I \subset k[x_1,\dots,x_N]$ be the ideal containing $(f_1,\dots,f_r)+(h)$ such that $\Spec k[x_1,\dots,x_N]/I$ is the nonsmooth locus of $V(h)$.@\footnote{See {\sc Singular}'s procedure \texttt{slocus} \cite[D.6.15.12]{GM20}. The procedure \texttt{slocus} computes the nonsmooth loci of affine varieties. If the ground field $k$ is perfect, singular loci and nonsmooth loci coincide.}@
$\text{\bf{if}}$ $1 \in I$ $\text{\bf{do}}$
  The nonsmooth locus of $V(h)$ is empty, hence $V(h)$ is a smooth hypersurface of $Y$.
  $\text{\bf{return}}$ $\T$
$\text{\bf{else}}$ 
  The nonsmooth locus of $V(h)$ is nonempty, hence $V(h)$ is a not smooth.
  $\text{\bf{return}}$ $\F$
\end{lstlisting}

\subsection{Examples and non-examples of smooth hypersurfaces}
The equation $y^3-x^2$ cuts a nonsmooth hypersurface out of $\mathbb{A}^2_k=\Spec k[x,y]$. The equation $z$ cuts a smooth hypersurface out of $\mathbb{A}^2_k \sqcup \mathbb{A}^2_k=\Spec k[x,y,z]/(x(x-1))$, but the equation $xz$ does not, because the subvariety cut out by $xz$ contains the component $\Spec k[x,y,z]/(x)$.

\subsection{Algorithm for Orthogonal idempotents on a smooth variety}
\label{subsec: idempotents}
Suppose $Y=\Spec A$ is a variety over a field $k$. Let $Y_1,\dots,Y_d$ be the connected components of $Y$, and set $A_i = \mathcal{O}(Y_i)$, so that $Y= \sqcup_i Y_i$ and $A= A_1 \times \cdots A_d$. The orthogonal idempotents of $Y$ are the elements $e_1,\dots,e_d\in A$ such that the image of $e_i$ in $A_1 \times \cdots A_d$ is the $i$th standard basis vector $(\dots,0,1,0,\dots)$, where the $i$th entry is $1$ and the other entries are zero \cite[Exercise 2.19]{Har77}. Note that $Y_i = \Spec A[e_i^{-1}]$.

If $Y=\Spec A$ is a smooth variety, we can explicitly compute the orthogonal idempotents of $Y$. The following algorithm, based on \cite[Proposition 1.10]{AM69}, computes the orthogonal idempotents of a smooth affine variety.

\Suppressnumber
\begin{lstlisting}[caption={smoothVarietyOrthogonalIdempotents},label={orthogonal idempotents}]
$\text{\bf{Input:}}$ $f_1,\dots,f_r \in k[x_1,\dots,x_N]$ such that $\Spec A$ is a smooth variety over $k$, where $A = k[x_1,\dots,x_N]/(f_1,\dots,f_r)$.
$\text{\bf{Output:}}$ Polynomials $e_1,\dots,e_d \in k[x_1,\dots,x_N]$ such that the irreducible components of $\Spec A$ are $\Spec A[e_1^{-1}],\dots,\Spec A[e_d^{-1}]$ @\Reactivatenumber@ 
Let $\frak{p}_1,\dots,\frak{p}_d \subset k[x_1,\dots,x_N]$ be the associated primes of the ideal $(f_1,\dots,f_r)$.@\footnote{See {\sc Singular}'s primary decomposition library \cite{PDSL19}}@
Because $\Spec A$ is a smooth variety, the associated primes $\{\frak{p}_i\}_i$ are the minimal primes of $A$ and $\{ \Spec A/\frak{p}_i \}_i$ are the irreducible components of $\Spec A$.
$\text{\bf{for}}$ $1 \leq i \leq d$ $\text{\bf{do}}$
  Initialization
  $e_i=1 \in k[x_1,\dots,x_N]$
  $\text{\bf{for}}$ each $1 \leq j \leq d$ such that $j \neq i$ $\text{\bf{do}}$
    Observe that the prime ideals $\frak{p}_i$ and $\frak{p}_j$ are coprime in $A$, that is, $\frak{p}_i  + \frak{p}_j =(1) \subset A$, because the closed subschemes $\Spec A/\frak{p}_i$ and $\Spec A/\frak{p}_j$ are disjoint. Since $\frak{p}_i$ and $\frak{p}_j$, as ideals of $k[x_1,\dots,x_N]$, contain $(f_1,\dots,f_r)$, it follows that $\frak{p}_i  + \frak{p}_j = (1) \subset k[x_1,\dots,x_N]$.
    Thus we can find polynomials $f,g \in k[x_1,\dots,x_N]$ such that $f \in \frak{p}_j$, $g \in \frak{p}_i$, and $f+g=1$ in $k[x_1,\dots,x_N]$@\footnote{Use {\sc Singular}'s procedure \texttt{lift} \cite[5.1.75]{DGPS19} to compute the polynomials $f$ and $g$.}@
    $f$ is the identity in $A/\frak{p}_i$, while $f$ is zero in $\Spec A/\frak{p}_j$.
    $e_i = e_i \cdot f\in k[x_1,\dots,x_N]$
  $e_i$ is the identity in $A/\frak{p}_i$, while $e_i$ is zero in $\Spec A/\frak{p}_j$ for each $j \neq i$. 
  Thus $\Spec A[e_i^{-1}] = \Spec A/\frak{p}_i$.
$\text{\bf{return}}$ $e_1,\dots,e_d$
\end{lstlisting}

\subsection{Example of orthogonal idempotents of smooth varieties}
The orthogonal idempotents of the disjoint union of the following three planes
$$
\Spec \tfrac{k[x,y,z]}{(x(x+1)(x+2))}=  \Spec\tfrac{k[x,y,z]}{(x)} \sqcup \Spec\tfrac{k[x,y,z]}{(x+1)} \sqcup \Spec\tfrac{k[x,y,z]}{(x+2)} 
$$
are the three polynomials $\frac{1}{2}(x+1)(x+2), \frac{1}{2}x(x+2), -x(x+1)$.



\begin{thebibliography}{ABMMOG14}

\bibitem[AFL20]{AFL20}
Dan Abramovich, Anne Fr\"{u}hbis-Kr\"{u}ger, Jonghyun Lee,
{\tt resweighted.lib}. {A} {\sc Singular} library for computing weighted resolution of singularities. To appear.

\bibitem[ATW20]{ATW20}
Dan Abramovich, Michael Temkin, and Jaros{\l}aw W{\l}odarczyk,
\emph{Functorial embedded resolution via weighted blowings up},
arXiv e-prints (2020), arXiv:1906.07106.

\bibitem[AM69]{AM69}
Michael F. Atiyah and Ian G. MacDonald,
\emph{Introduction to Commutative Algebra},
Addison-Wesley, 
London (1969).

\bibitem[BM97]{BM97}
Edward Bierstone and Pierre Milman, 
\emph{Canonical desingularization in characteristic zero by blowing-up the maxima strata of a local invariant}, 
Invent. Math. 128 (1997), 207?302.

\bibitem[BS00a]{BS00a}
G{\'a}bor Bodn{\'a}r and Josef Schicho, 
\emph{Automated resolution of singularities for hypersurfaces}, 
J. Symbolic Comput. 30 (2000), no. 4, 401428.

\bibitem[BS00b]{BS00b}
G{\'a}bor Bodn{\'a}r and Josef Schicho, 
\emph{A computer program for the resolution of singularities}, 
Resolution of singularities (Obergurgl, 1997), Progr. Math., vol. 181, Birkh{\"a}user, Basel, 2000, pp. 231-238.

\bibitem[BDFPRR18]{BDFPRR18}
Janko Boehm, Wolfram Decker, Anne Fr\"{u}hbis-Kr\"{u}ger, Franz-Josef Pfreundt, Mirko Rahn, and Lukas Ristau,
\emph{Towards Massively Parallel Computations in Algebraic Geometry},
arXiv e-prints (2018), arXiv:1808.09727.

\bibitem[BEV05]{BEV05}
Ana Bravo, Santiago Encinas, Orlando Villamayor,  
\emph{A Simplified Proof of Desingularisation and Applications},
Rev. Math. Iberoamericana 21, 349?458 (2005).

\bibitem[BV01]{BV01}
Ana Bravo and Orlando Villamayor, 
\emph{Strengthening the theorem of embedded desingularization}, 
Math. Res. Lett. 8 (2001), no. 1-2, 79-89.

\bibitem[CLO15]{CLO15}
David A. Cox, John Little, and Donal O?Shea, 
\emph{Ideals, Varieties, and Algorithms}, 
Fourth Edition, Undergraduates Texts In Mathematics, Springer, New York, 2015.
MR 3330490

\bibitem[DGPS19]{DGPS19}
Wolfram Decker, Gert-Martin Greuel, Gerhard Pfister, Hans Sch{\"o}nemann,    
\newblock {\sc Singular} {4-1-2} --- {A} computer algebra system for polynomial computations.
\newblock {http://www.singular.uni-kl.de} (2020).

\bibitem[EV98]{EV98}
Santiago Encinas and Orlando Villamayor, 
\emph{Good points and constructive resolution of singularities},
Acta Math. 181 (1998), 109-158.

\bibitem[F07]{F07}
Anne Fr\"{u}hbis-Kr\"{u}ger,
\emph{Computational Aspects of Singularities}, 
in J.-P. Brasselet, J.Damon et al.: Singularities in Geometry and Topology, World Scientific Publishing (2007), 253?327.


\bibitem[FP19]{FP19}
Anne Fr\"{u}hbis-Kr\"{u}ger and Gerhard Pfister,
{\tt resolve.lib}. {A} {\sc Singular} {4-1-2} library for computing resolution of singularities (2019). 

\bibitem[GM20]{GM20}
Gert-Martin Greuel and Bernd Martin,
{\tt sing.lib}. {A} {\sc Singular} {4-1-2} library for computing invariants of singularities (2020). 

\bibitem[EGAIV]{EGAIV}
Alexander Grothendieck and Jean Dieudonn{\'e},
\emph{\'{E}lements de g\'{e}om\'{e}trie alg\'{e}brique IV: \'{E}tude locale des sch\'{e}mas et des morphismes de sch\'{e}mas, Quatri\`{e}me partie},
Publications Math\'{e}matiques de L' I.H.\'{E}.S., 32 (1967).

\bibitem[Har77]{Har77}
Robin Hartshorne,
\emph{Algebraic Geometry},
Springer-Verlag, New York, 1977, 
Graduate Texts in Mathematics, No. 52. 
MR 0463157 (57 \#3116)

\bibitem[Hir64]{Hir64}
Heisuke Hironaka, 
\emph{Resolution of singularities of an algebraic variety over a field of characteristic zero}.
I, II, Ann. Of Math. (2) 79 (1964), 109-203; 
ibid. (2) 79 (1964), 205-326. 
MR 0199184

\bibitem[Kol07]{Kol07}
J{\'a}nos Koll{\'a}r, 
\emph{Lectures on resolution of singularities}, 
Annals of Mathematics Studies, vol. 166, Princeton University Press, Princeton, NJ, 2007.

\bibitem[Liu06]{Liu06}
Qing Liu,
\emph{Algebraic Geometry and Arithmetic Curves},
Oxford Graduate Texts in Mathematics. Oxford University Press, Oxford, 2006.

\bibitem[Mat89]{Mat89}
Hideyuki Matsumura,
\emph{Commutative ring theory},
Cambridge Studies in Advanced Mathematics, vol. 8, Cambridge University Press, Cambridge, 1989.

\bibitem[McQ19]{McQ19}
Michael McQuillan, 
\emph{Very fast, very functorial, and very easy resolution of singularities}, 
2019, preprint, written with the collaboration of Gianluca Marzo.

\bibitem[PDSL19]{PDSL19}
Gerhard Pfister, Wolfram Decker, Hans Schoenemann, and Santiago Laplagne,
{\tt primdec.lib}. {A} {\sc Singular} {4-1-2} library for computing primary decomposition (2019). 

\bibitem[Que20]{Que20}
Ming Hao Quek,
\emph{Logarithmic resolution via weighted toroidal blowings up},
arXiv e-prints (2020), arXiv:1906.07106.

\bibitem[SH18]{SH18}
Irene Swanson and Craig Huneke,
\emph{Integral Closure of Ideals, Rings, and Modules},
online version \url{http://people.reed.edu/~iswanson/book/SwansonHuneke.pdf} February 2018.

\bibitem[Stacks20]{Stacks20}
{The {Stacks project authors}},
{The Stacks project},
{\url{https://stacks.math.columbia.edu}},
{2020}.

\bibitem[Vak17]{Vak17}
Ravi Vakil,
\emph{The Rising Sea: Foundations of Algebraic Geometry},
November 18, 2017 draft,
\url{http://math.stanford.edu/~vakil/216blog/FOAGnov1817public.pdf}.

\bibitem[Vil89]{Vil89}
Orlando Villamayor, 
\emph{Constructiveness of Hironaka's resolution}, 
Ann. Sci. {\'E}cole Norm. Sup. (4) 22 (1989), no. 1, 1-32.

\bibitem[Vil92]{Vil92}
Orlando Villamayor, 
\emph{Patching local uniformizations}, 
Ann. Sci. {\'E}cole Norm. Sup. (4) 25 (1992), no. 6, 629677.

\bibitem[W{\l}o05]{Wlo05}
Jaros{\l}aw W{\l}odarczyk,
\emph{Simple Hironaka resolution in characteristic zero}, 
J. Amer. Math. Soc. 18 (2005), no. 4, 779?822 (electronic). MR 2163383

\end{thebibliography}
\end{document}